\documentclass[12pt,twoside]{article}

\usepackage{amsfonts}
\usepackage{amssymb}
\usepackage{amsmath}
\usepackage{graphicx}

\newcommand{\ignore}[1]{}

\def\@begintheorem#1#2{\par\bgroup{\sc #1\ #2. }\it\ignorespaces}
\def\@opargbegintheorem#1#2#3{\par\bgroup{\sc #1\ #2\ (#3). } \it\ignorespaces}
\def\@endtheorem{\egroup}
\newtheorem{theorem}{Theorem}[section]
\newtheorem{corollary}[theorem]{Corollary}
\newtheorem{lemma}[theorem]{Lemma}
\newtheorem{proposition}[theorem]{Proposition}
\newtheorem{example}[theorem]{Example}
\newtheorem{algorithm}[theorem]{Algorithm}
\newtheorem{definition}[theorem]{Definition}
\newcommand{\bt}[1]{\begin{theorem}\label{#1}}
\newcommand{\bc}[1]{\begin{corollary}\label{#1}}
\newcommand{\bl}[1]{\begin{lemma}\label{#1}}
\newcommand{\bp}[1]{\begin{proposition}\label{#1}}
\newcommand{\be}[1]{\begin{example}\label{#1}}
\newcommand{\ba}[1]{\begin{algorithm}\rm\label{#1}}
\newcommand{\bd}[1]{\begin{definition}\rm\label{#1}}
\newcommand{\bpr}{\noindent {\em Proof. }}
\newcommand{\et}{\end{theorem}}
\newcommand{\ec}{\end{corollary}}
\newcommand{\el}{\end{lemma}}
\newcommand{\ep}{\end{proposition}}
\newcommand{\ee}{\end{example}}
\newcommand{\ea}{\end{algorithm}}
\newcommand{\ed}{\end{definition}}
\newcommand{\epr}{{\ \vbox{\hrule\hbox{%
\vrule height1.3ex\hskip0.8ex\vrule}\hrule}}\\\par}
\newcommand{\mepr}{{\ \ \ \vbox{\hrule\hbox{%
\vrule height1.3ex\hskip0.8ex\vrule}\hrule}}}

\def\R{\mathbb{R}}
\def\Z{\mathbb{Z}}

\def \G {{\cal G}}
\def \l {\langle}
\def \r {\rangle}

\def \L {{\cal L}}
\def \G {{\cal G}}

\def \l {\langle}
\def \r {\rangle}

\def \A {A^{(n)}}

\def \conv {{\rm conv}}

\def \type {{\rm type}}

\begin{document}

\title{\bf Theory and Applications of \\
N-Fold Integer Programming}

\author{
Shmuel Onn
\thanks{Supported in part by a grant from ISF - the Israel Science Foundation}
}

\date{}
\maketitle

\begin{abstract}
We overview our recently introduced theory of $n$-fold integer programming
which enables the polynomial time solution of fundamental linear and nonlinear
integer programming problems in variable dimension. We demonstrate its
power by obtaining the first polynomial time algorithms in several
application areas including multicommodity flows and
privacy in statistical databases.
\end{abstract}

\section{Introduction}
\label{i}

Linear integer programming is the following fundamental optimization problem,
$$\min\,\left\{wx\ :\ x\in\Z^n\,,\ Ax=b\,,\ l\leq x\leq u\right\}\ ,$$
where $A$ is an integer $m\times n$ matrix, $b\in\Z^m$, and $l,u\in\Z_{\infty}^n$
with $\Z_{\infty}:=\Z\uplus\{\pm\infty\}$.
It is generally NP-hard, but polynomial time solvable in two
fundamental situations: the dimension is fixed \cite{Len};
the underlying matrix is totally unimodular \cite{HK}.

Recently, in \cite{DHOW}, a new fundamental polynomial time solvable
situation was discovered. We proceed to describe this class
of so-termed {\em $n$-fold integer programs}.

An {\em $(r,s)\times t$ bimatrix}
is a matrix $A$ consisting of two blocks $A_1$, $A_2$, with $A_1$
its $r\times t$ submatrix consisting of the first $r$ rows and $A_2$
its $s\times t$ submatrix consisting of the last $s$ rows. The
{\em $n$-fold product} of $A$ is the following $(r+ns)\times nt$ matrix,
$$A^{(n)}\quad:=\quad
\left(
\begin{array}{cccc}
  A_1    & A_1    & \cdots & A_1    \\
  A_2    & 0      & \cdots & 0      \\
  0      & A_2    & \cdots & 0      \\
  \vdots & \vdots & \ddots & \vdots \\
  0      & 0      & \cdots & A_2    \\
\end{array}
\right)\quad .
$$
The following result of \cite{DHOW} asserts
that $n$-fold integer programs are efficiently solvable.
\bt{NFold1}{\bf \cite{DHOW}}
For each fixed integer $(r,s)\times t$ bimatrix $A$, there is
an algorithm that, given positive integer $n$, bounds $l,u\in\Z_{\infty}^{nt}$,
$b\in\Z^{r+ns}$, and $w\in\Z^{nt}$, solves in time which is polynomial in $n$
and in the binary-encoding length $\l l,u,b,w \r$ of the rest of the data,
the following so-termed linear $n$-fold integer programming problem,
$$\min\,\left\{wx\ :\ x\in\Z^{nt}\,,\ A^{(n)} x=b\,,\ l\leq x\leq u\right\}\ .$$
\et

\vskip.25cm
Some explanatory notes are in order. First, the dimension
of an $n$-fold integer program is $nt$ and is variable.
Second, $n$-fold products $A^{(n)}$ are highly non totally unimodular:
the $n$-fold product of the simple $(0,1)\times 1$ bimatrix
with $A_1$ empty and $A_2:=2$ satisfies $A^{(n)}=2 I_n$ and has exponential
determinant $2^n$. So this is indeed a class of programs which cannot
be solved by methods of fixed dimension or totally unimodular matrices.
Third, this class of programs turns out to be very natural and has
numerous applications, the most generic being to integer optimization
over multidimensional tables (see \S \ref{a}). In fact it is {\em universal}:
the results of \cite{DO2} imply that {\em every} integer program {\em is}
an $n$-fold program over some simple bimatrix $A$ (see \S \ref{d}).

The above theorem extends to $n$-fold integer programming with
nonlinear objective functions as well. The following results,
from \cite{HOW1}, \cite{DHORW} and \cite{HOW2}, assert that the
minimization and maximization of broad classes of convex functions over
$n$-fold integer programs can also be done in polynomial time.
The function $f$ is presented either by a {\em comparison oracle} that
for any two vectors $x,y$ can answer whether or not $f(x)\leq f(y)$, or by an
{\em evaluation oracle} that for any vector $x$ can return $f(x)$.

In the next theorem, $f$ is {\em separable convex}, namely $f(x)=\sum_i f_i(x_i)$
with each $f_i$ univariate convex. Like linear forms, such functions
can be minimized over totally unimodular programs \cite{HSh}.
We show that they can also be efficiently minimized over $n$-fold programs.
The running time depends also on $\log {\hat f}$ with ${\hat f}$ the maximum
value of $|f(x)|$ over the feasible set (which need not be part of the input).
\bt{NFold2}{\bf \cite{HOW1}}
For each fixed integer $(r,s)\times t$ bimatrix $A$, there is an algorithm
that, given $n$, $l,u\in\Z_{\infty}^{nt}$, $b\in\Z^{r+ns}$,
and separable convex $f:\Z^{nt}\rightarrow\Z$
presented by a comparison oracle, solves in time polynomial in $n$ and
$\l l,u,b,{\hat f}\r$, the program
\begin{equation*}
\min\left\{f(x)\ :\ x\in\Z^{nt}\,,\ A^{(n)} x=b\,,\ l\leq x\leq u\right\}\ .
\end{equation*}
\et

\vskip.25cm
An important natural special case of Theorem \ref{NFold2} is the following result that
concerns finding a feasible point which is $l_p$-closest to a given desired goal point.
\bt{NFold3}{\bf \cite{HOW1}}
For each fixed integer $(r,s)\times t$ bimatrix $A$, there is an algorithm
that, given positive integers $n$ and $p$, $l,u\in\Z_{\infty}^{nt}$, $b\in\Z^{r+ns}$,
and ${\hat x}\in\Z^{nt}$, solves in time polynomial in $n$, $p$,
and $\l l,u,b,{\hat x}\r$, the following distance minimization program,
\begin{equation}\label{distance_equation}
\min\,\{\|x-{\hat x}\|_p\ :\ x\in\Z^{nt},\ \A x=b,\ l\leq x\leq u\}\ .
\end{equation}
For $p=\infty$ the problem (\ref{distance_equation})
can be solved in time polynomial in $n$ and $\l l,u,b,{\hat x}\r$.
\et

\vskip.25cm
The next result concerns the {\em maximization} of a convex function
of the composite form $f(Wx)$, with $f:\Z^d\rightarrow\Z$ convex
and $W$ an integer matrix with $d$ rows.
\bt{NFold4}{\bf \cite{DHORW}}
For each fixed $d$ and $(r,s)\times t$ integer bimatrix $A$,
there is an algorithm that, given $n$, bounds
$l,u\in\Z_{\infty}^{nt}$, integer $d\times {nt}$ matrix $W$,
$b\in\Z^{r+ns}$, and convex function $f:\Z^d\rightarrow\R$ presented
by a comparison oracle, solves in time polynomial in $n$ and $\l W,l,u,b \r$,
the convex $n$-fold integer maximization program
$$\max\{f(Wx)\ :\ x\in\Z^{nt}\,,\ A^{(n)}x=b\,,\ l\leq x\leq u\}\ .$$
\et

\vskip.25cm
Finally, we have the following broad extension of Theorem \ref{NFold2}
where the objective can include a composite term $f(Wx)$,
with $f:\Z^d\rightarrow\Z$ separable convex and $W$
an integer matrix with $d$ rows, and where also inequalities on $Wx$
can be included. As before, ${\hat f},{\hat g}$ denote the maximum
values of $|f(Wx)|,|g(x)|$ over the feasible set.
\bt{NFold5}{\bf \cite{HOW2}}
For each fixed integer $(r,s)\times t$ bimatrix $A$ and integer $(p,q)\times t$
bimatrix $W$, there is an algorithm that, given $n$,
$l,u\in\Z_{\infty}^{nt}$, ${\hat l},{\hat u}\in\Z_{\infty}^{p+nq}$, $b\in\Z^{r+ns}$,
and separable convex
functions $f:\Z^{p+nq}\rightarrow\Z$, $g:\Z^{nt}\rightarrow\Z$ presented
by evaluation oracles, solves in time
polynomial in $n$ and $\l l,u,{\hat l},{\hat u},b,{\hat f},{\hat g}\r$,
the generalized program
\begin{equation*}\label{NFoldStrongSeparableConvexEquation}
\min\left\{f(W^{(n)}x)+g(x)\ :\ x\in\Z^{nt}\,,\ A^{(n)}x=b\,,\
{\hat l}\leq W^{(n)}x\leq{\hat u}\,,\ l\leq x\leq u\right\}\ .
\end{equation*}
\et

\vskip0.5cm
The article is organized as follows.
In Section 2 we discuss some of the many applications of this
theory and use Theorems \ref{NFold1}--\ref{NFold5} to obtain the
first polynomial time algorithms for these applications.
In Section 3 we provide a concise development of the theory
of $n$-fold integer programming and prove our
Theorems \ref{NFold1}--\ref{NFold5}. Sections 2 and 3
can be read in any order. We conclude in Section 4 with a
discussion of the universality of $n$-fold integer programming
and of a new (di)-graph invariant, about which very little is known,
that is important in understanding the complexity of our algorithms.
Further discussion of $n$-fold integer programming within the broader context of
{\em nonlinear discrete optimization} can be found in \cite{Onn2} and \cite{Onn3}.

\section{Applications}
\label{a}

\subsection{Multiway Tables}
\label{mt}

Multiway tables occur naturally in any context
involving multiply-indexed variables. They have been studied extensively
in mathematical programming in the context of high dimensional transportation
problems (see \cite{Vla,YKK} and the references therein) and in statistics
in the context of disclosure control and privacy in statistical databases
(see \cite{Cox,FR} and the references therein).
The theory of $n$-fold integer programming provides the first polynomial time
algorithms for multiway table problems in these two contexts,
which are discussed in \S \ref{mitp} and \S \ref{pisd} respectively.

\begin{figure}[hbt]
\hskip-1.1cm
\includegraphics[scale=0.58]{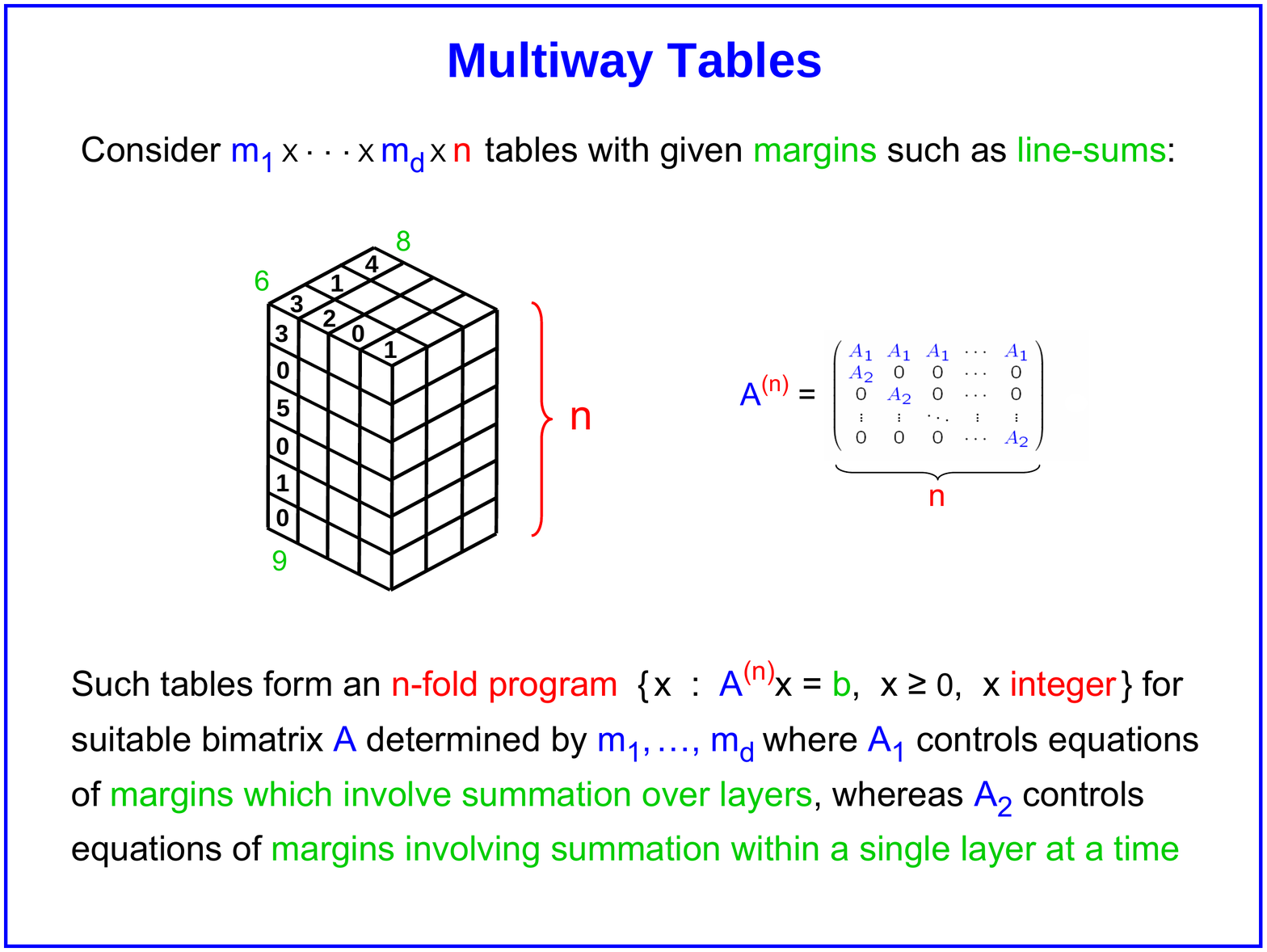}
\caption{Multiway Tables}
\label{multiway-tables-figure}
\end{figure}

We start with some terminology and background that will be used in
the sequel. A {\em $d$-way table} is an $m_1\times\cdots\times m_d$
array $x=(x_{i_1,\dots,i_d})$ of nonnegative integers.
A {\em $d$-way transportation polytope} ({\em $d$-way polytope}
for brevity) is the set of $m_1\times\cdots\times m_d$
nonnegative arrays $x=(x_{i_1,\dots,i_d})$ with specified
sums of entries over some of their lower dimensional subarrays
({\em margins} in statistics). The $d$-way tables with specified
margins are the integer points in the $d$-way polytope.
For example (see Figure \ref{multiway-tables-figure}),
the $3$-way polytope of $l\times m\times n$ arrays
with specified line-sums ($2$-margins) is
\begin{equation*}
T\ \ :=\ \ \left\{x\in\R_+^{l\times  m\times n}\ :\ \sum_i x_{i,j,k}=v_{*,j,k}
\,,\ \sum_j x_{i,j,k}=v_{i,*,k}\,,\ \sum_k x_{i,j,k}=v_{i,j,*}\right\}\ \ ,
\end{equation*}
where the specified line-sums are $mn+ln+lm$ given nonnegative integer numbers
$$v_{*,j,k}\,,\ \ v_{i,*,k}\,,\ \ v_{i,j,*}\,,\quad\quad
1\leq i\leq l\,,\ \ 1\leq j\leq m\,,\ \ 1\leq k\leq n\quad .$$
Our results hold for $k$-margins
for any $0\leq k\leq d$, and much more generally for any so-called
{\em hierarchical family} of margins. For simplicity of the exposition, however,
we restrict attention here to line-sums, that is, $(d-1)$-margins, only.

We conclude this preparation with the {\em universality theorem}
for multiway tables and polytopes. It provides a powerful
tool in establishing the presumable limits of polynomial time
solvability of table problems, and will be used in \S \ref{mitp} and \S \ref{pisd}
to contrast the polynomial time solvability attainable by $n$-fold integer programming.
\bt{TableUniversality}{\bf \cite{DO2}}
Every rational polytope $P=\{y\in\R_+^d\,:\,Ay=b\}$ is in polynomial
time computable integer preserving bijection with some
$l\times m\times 3$ line-sum polytope
\begin{equation*}
T\ \ =\ \ \left\{x\in\R_+^{l\times  m\times 3}\ :\ \sum_i x_{i,j,k}=v_{*,j,k}
\,,\ \sum_j x_{i,j,k}=v_{i,*,k}\,,\ \sum_k x_{i,j,k}=v_{i,j,*}\right\}\ \ .
\end{equation*}
\et

\subsubsection{Multi-index transportation problems}
\label{mitp}

The {\em multi-index transportation problem} of Motzkin \cite{Mot} is the integer
programming problem over multiway tables with specified margins.
For line-sums it is the program
\begin{equation*}\label{line_sum_problem}
\min\left\{wx\ :\ x\in\Z_+^{m_1\times\cdots \times m_d}\ :\
\sum_{i_1} x_{i_1,\dots,i_d}=v_{*,i_2,\dots, i_d}\,,
\ \dots\,,\ \sum_{i_d} x_{i_1,\dots,i_d}=v_{i_1,\dots,i_{d-1},*}\right\}.
\end{equation*}
For $d=2$ this program is totally unimodular and can be solved in
polynomial time. However, already for $d=3$ it is generally not,
and the problem is much harder. Consider the problem
over $l\times m\times n$ tables. If $l,m,n$ are all fixed
then the problem is solvable in polynomial time (in the natural
binary-encoding length of the line-sums), but even in this very restricted
situation one needs off-hand the algorithm of integer programming in
fixed dimension $lmn$. If $l,m,n$ are all variable then the problem
is NP-hard \cite{IJ}. The in-between cases are much more delicate and
were resolved only recently. If two sides are variable and one is fixed
then the problem is still NP-hard \cite{DO1}; moreover,
Theorem \ref{TableUniversality} implies that it is NP-hard even over
$l\times m\times 3$ tables with fixed $n=3$.
Finally, if two sides are fixed and one is variable,
then the problem can be solved in polynomial time by $n$-fold
integer programming. Note that even over $3\times 3\times n$ tables,
the only solution of the problem available to-date
is the one given below using $n$-fold integer programming.

The polynomial time solvability of the multi-index transportation
problem when one side is variable and the others are fixed extends to
any dimension $d$. We have the following important result on
the multi-index transportation problem.
\bt{Multiway}{\bf \cite{DHOW}}
For every fixed $d,m_1,\dots,m_d$, there is an algorithm that, given $n$,
integer $m_1\times\cdots \times m_d\times n$ cost $w$, and integer line-sums
$v=((v_{*,i_2,\dots, i_{d+1}}),\dots,(v_{i_1,\dots,i_d,*}))$, solves in time
polynomial in $n$ and $\l w,v \r$, the $(d+1)$-index transportation problem
$$\min\left\{wx:x\in\Z_+^{m_1\times\cdots\times m_d\times n}:
\sum_{i_1} x_{i_1,\dots,i_{d+1}}=v_{*,i_2,\dots, i_{d+1}},\dots,
\sum_{i_{d+1}} x_{i_1,\dots,i_{d+1}}=v_{i_1,\dots,i_d,*}\right\}.$$
\et
\bpr
Re-index arrays as $x=(x^1,\dots,x^n)$ with each
$x^{i_{d+1}}=(x_{i_1,\dots,i_d,i_{d+1}})$ a suitably indexed $m_1m_2\cdots m_d$ vector
representing the $i_{d+1}$-th layer of $x$. Similarly re-index the array $w$.
Let $t:=r:=m_1m_2\cdots m_d$ and
$s:=n\left(m_2\cdots m_d\,+\,\cdots\,+\,m_1\cdots m_{d-1}\right)$.
Let $b:=(b^0,b^1,\dots,b^n)\in\Z^{r+ns}$, where
$b^0:=(v_{i_1,\dots,i_d,*})$ and for $i_{d+1}=1,\dots,n$,
$$b^{i_{d+1}}\ :=\ \left((v_{*,i_2,\dots,i_d,i_{d+1}}),\dots,
(v_{i_1,\dots,i_{d-1},*,i_{d+1}})\right)\ .$$
Let $A$ be the $(t,s)\times t$ bimatrix with first block $A_1:=I_t$
the $t\times t$ identity matrix and second block $A_2$ a matrix
defining the line-sum equations on $m_1\times\cdots\times m_d$ arrays.
Then the equations $A_1(\sum_{i_{d+1}} x^{i_{d+1}})=b^0$ represent
the line-sum equations $\sum_{i_{d+1}} x_{i_1,\dots,i_{d+1}}=v_{i_1,\dots,i_d,*}$
where summations over layers occur, whereas the equations $A_2x^{i_{d+1}}=b^{i_{d+1}}$
for $i_{d+1}=1,\dots,n$ represent all other line-sum equations,
where summations are within a single layer at a time. Therefore
the multi-index transportation problem
is encoded as the $n$-fold integer programming problem
$$\min\,\{wx\ :\ x\in\Z^{nt},\ A^{(n)}x=b,\ x\geq 0\}\ .$$
Using the algorithm of Theorem \ref{NFold1},
this $n$-fold integer program, and hence the given multi-index
transportation problem, can be solved in polynomial time.
\epr
This proof extends immediately to multi-index transportation problems
with nonlinear objective functions of the forms in
Theorems \ref{NFold2}--\ref{NFold5}. Moreover, as mentioned before,
a similar proof shows that multi-index transportation problems with $k$-margin
constraints, and more generally, hierarchical margin constraints, can be encoded
as $n$-fold integer programming problems as well. We state this as a corollary.
\bc{CorollaryTransportation}{\bf \cite{DHORW}}
For every fixed $d$ and $m_1,\dots,m_d$, the nonlinear multi-index transportation
problem, with any hierarchical margin constraints, over $(d+1)$-way tables of
format $m_1\times\cdots\times m_d\times n$ with variable $n$ layers, are
polynomial time solvable.
\ec

\subsubsection{Privacy in statistical databases}
\label{pisd}

A common practice in the disclosure of sensitive data contained in a
multiway table is to release some of the table margins rather than the
entries of the table. Once the margins are released, the security of any specific
entry of the table is related to the set of possible values that can
occur in that entry in all tables having the same margins as those
of the source table in the database. In particular, if this set consists
of a unique value, that of the source table, then this entry can be exposed
and privacy can be violated. This raises the following fundamental problem.

\vskip.2cm\noindent{\bf Entry uniqueness problem:} Given
hierarchical margin family and entry index, is the value which
can occur in that entry in all tables with these margins, unique?

\vskip.2cm\noindent
The complexity of this problem turns out to behave in analogy to the complexity of
the multi-index transportation problem discussed in \S \ref{mitp}. Consider the problem
for $d=3$ over $l\times m\times n$ tables. It is polynomial time decidable when
$l,m,n$ are all fixed, and coNP-complete when $l,m,n$ are all variable \cite{IJ}.
We discuss next in more detail the in-between cases
which are more delicate and were settled only recently.

\begin{figure}[hbt]
\hskip-1.1cm
\includegraphics[scale=0.58]{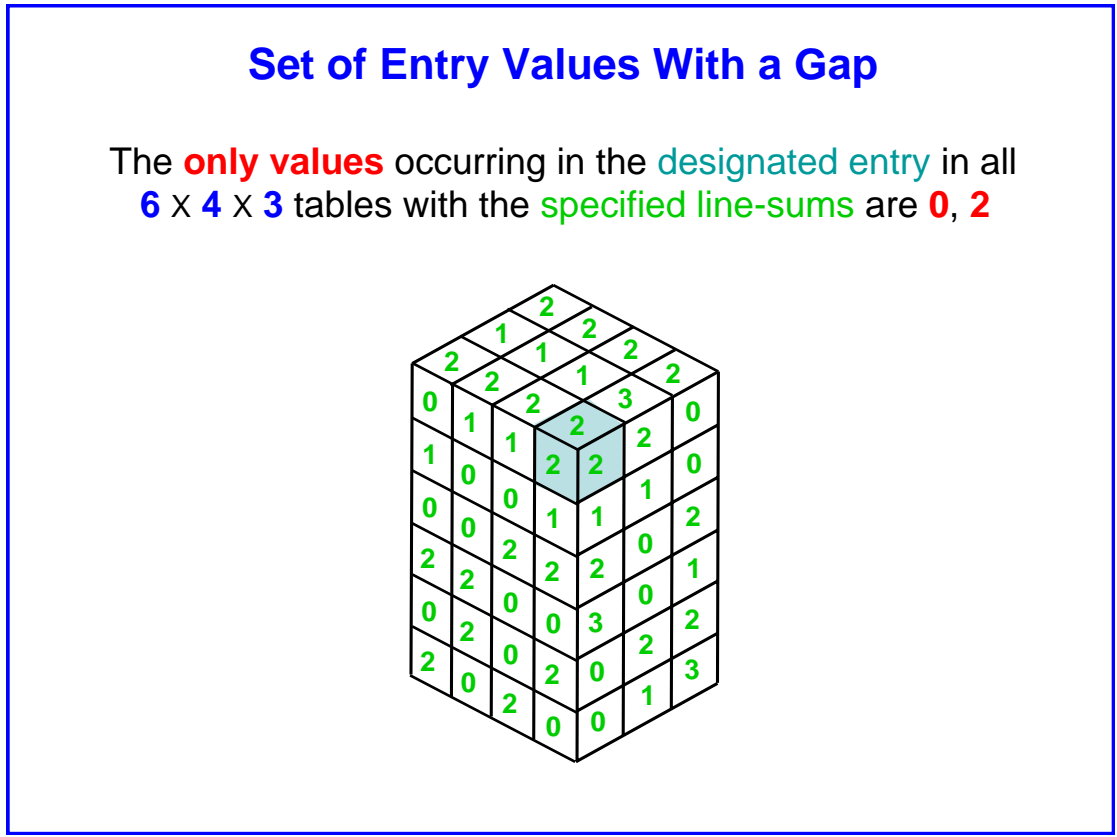}
\caption{Set of Entry Values With a Gap}
\label{Set-of-Entry-Values-With-a-Gap}
\end{figure}

If two sides are variable and one is fixed then the problem is still
coNP-complete, even over $l\times m\times 3$ tables with fixed $n=3$
\cite{Onn1}. Moreover, Theorem \ref{TableUniversality}
implies that {\em any set of nonnegative integers} is the set of values
of an entry of some $l\times m\times 3$ tables with some specified line-sums.
Figure \ref{Set-of-Entry-Values-With-a-Gap} gives an example of
line-sums for $6\times 4\times 3$ tables where one entry attains
the set of values $\{0,2\}$ which has a {\em gap}.

\bt{EntryUniversality}{\bf \cite{DO3}}
For every finite set $S\subset\Z_+$ of nonnegative integers, there
exist $l,m$, and line-sums for $l\times m\times 3$ tables, such
that the set of values that occur in some fixed entry in all
$l\times m\times 3$ tables that have these line-sums, is precisely $S$.
\et
\bpr
Consider any finite set $S=\{s_1,\dots,s_h\}\subset\Z_+$. Consider the polytope
$$P\ \ :=\ \ \{y\in\R_+^{h+1}\ :\
y_0-\sum_{j=1}^h s_jy_j=0\,,\ \sum_{j=1}^h y_j=1\,\}\ \ .$$
By Theorem \ref{TableUniversality}, there are $l,m$,
and $l\times m\times 3$ polytope $T$ with line-sums
$$v_{*,j,k}\,,\ \ v_{i,*,k}\,,\ \ v_{i,j,*}\,,\quad\quad
1\leq i\leq l\,,\ \ 1\leq j\leq m\,,\ \ 1\leq k\leq 3\quad\ , $$
such that the integer points in $T$, which are precisely the
$l\times m\times 3$ tables with these line-sums, are in bijection
with the integer points in $P$. Moreover (see \cite{DO2}),
this bijection is obtained by a simple projection from
$\R^{l\times m\times 3}$ to $\R^{h+1}$ that erases all
but some $h+1$ coordinates. Let $x_{i,j,k}$ be the coordinate
that is mapped to $y_0$. Then the set of values that this entry
attains in all tables with these line-sums is, as desired,
\begin{eqnarray*}
\left\{x_{i,j,k}\, :\, x\in T\cap\Z^{l\times m\times 3}\right\}
\ =\ \left\{y_0\, :\, y\in P\cap\Z^{h+1}\right\}\ \ =\ \ S\ . \ \
\mepr
\end{eqnarray*}

\vskip.4cm
Finally, if two sides are fixed and one is variable, then entry uniqueness
can be decided in polynomial time by $n$-fold integer programming.
Note that even over $3\times 3\times n$ tables, the only solution of the
problem available to-date is the one below.

The polynomial time decidability of the problem when one side
is variable and the others are fixed extends to any dimension $d$.
It also extends to any hierarchical family of margins,
but for simplicity we state it only for line-sums, as follows.

\bt{EntryUniqueness}{\bf \cite{Onn1}}
For every fixed $d,m_1,\dots,m_d$, there is an algorithm that, given $n$,
integer line-sums $v=((v_{*,i_2,\dots, i_{d+1}}),\dots,(v_{i_1,\dots,i_d,*}))$,
and entry index $(k_1,\dots,k_{d+1})$, solves in time which is polynomial in $n$
and $\l v \r$, the corresponding entry uniqueness problem, of deciding if
the entry $x_{k_1,\dots,k_{d+1}}$ is the same in all $(d+1)$-tables in the set
$$S\, :=\, \left\{x\in\Z_+^{m_1\times\cdots\times m_d\times n}\ :\
\sum_{i_1} x_{i_1,\dots,i_{d+1}}=v_{*,i_2,\dots, i_{d+1}},\dots,
\sum_{i_{d+1}} x_{i_1,\dots,i_{d+1}}=v_{i_1,\dots,i_d,*}\right\}\, .$$
\et
\bpr
By Theorem \ref{Multiway} we can solve in polynomial time both $n$-fold programs
$$l\quad:=\quad\min\left\{x_{k_1,\dots,k_{d+1}}\ :\ x\in S\right\}\ ,$$
$$u\quad:=\quad\max\left\{x_{k_1,\dots,k_{d+1}}\ :\ x\in S\right\}\ .$$
Clearly, entry $x_{k_1,\dots,k_{d+1}}$ has the same value in all
tables with the given line-sums if and only if $l=u$,
which can therefore be tested in polynomial time.
\epr
The algorithm of Theorem \ref{EntryUniqueness} and its extension
to any family of hierarchical margins allow statistical agencies
to efficiently check possible margins before disclosure:
if an entry value is not unique then disclosure may be
assumed secure, whereas if the value is unique then
disclosure may be risky and fewer margins should be released.

We note that {\em long} tables, with one side much larger
than the others, often arise in practical applications.
For instance, in health statistical tables, the long factor may be
the age of an individual, whereas other factors may be binary
(yes-no) or ternary (subnormal, normal, and supnormal).
Moreover, it is always possible to merge categories of factors,
with the resulting coarser tables approximating the original ones,
making the algorithm of Theorem \ref{EntryUniqueness} applicable.

Finally, we describe a procedure based on
a suitable adaptation of the algorithm of Theorem \ref{EntryUniqueness},
that constructs the entire set of values that can occur in a
specified entry, rather than just decides its uniqueness. Here $S$
is the set of tables satisfying the given (hierarchical) margins,
and the running time is output-sensitive, that is, polynomial in
the input encoding plus the number of elements in the output set.

\vskip.2cm\noindent{\bf Procedure for constructing the set of values in an entry:}
\begin{enumerate}
\item
Initialize $l:=-\infty$, $u:=\infty$, and $E:=\emptyset$.
\item
Solve in polynomial time the following linear $n$-fold integer programs:
$${\hat l}\ :=\ \min\left\{x_{k_1,\dots,k_{d+1}}\, :\,
l\leq x_{k_1,\dots,k_{d+1}}\leq u\,,\quad x\in S\right\}\ ,$$
$${\hat u}\ :=\ \max\left\{x_{k_1,\dots,k_{d+1}}\, :\,
l\leq x_{k_1,\dots,k_{d+1}}\leq u\,,\quad x\in S\right\}\ .$$
\item
If the problems in Step 2 are feasible then update
$l:={\hat l}+1$, $u:={\hat u}-1$, $E:=E\uplus \{{\hat l},{\hat u}\}$,
and repeat Step 2, else stop and output the set of values $E$.
\end{enumerate}

\subsection{Multicommodity Flows}
\label{mf}

\begin{figure}[hbt]
\hskip-1.1cm
\includegraphics[scale=0.58]{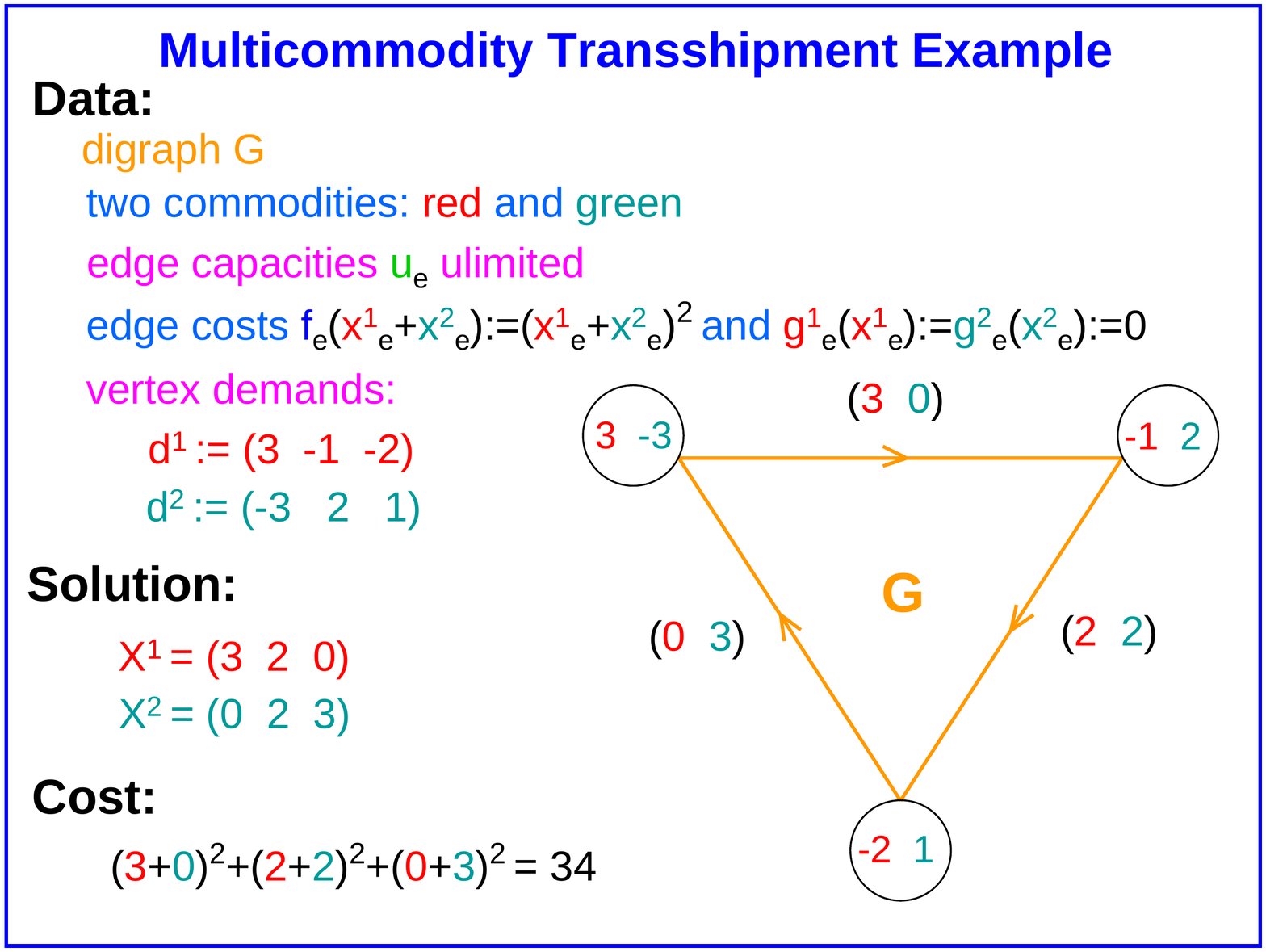}
\caption{Multicommodity Transshipment Example}
\label{multi-transshipment-figure}
\end{figure}

The multicommodity transshipment problem is a very general flow problem
which seeks minimum cost routing of several discrete commodities
over a digraph subject to vertex demand and edge capacity constraints.
The data for the problem is as follows
(see Figure \ref{multi-transshipment-figure} for a small example).
There is a digraph $G$ with $s$ vertices and $t$ edges.
There are $l$ types of commodities.
Each commodity has a demand vector $d^k\in\Z^s$ with $d^k_v$
the demand for commodity $k$ at vertex $v$ (interpreted as supply
when positive and consumption when negative).
Each edge $e$ has a capacity $u_e$
(upper bound on the combined flow of all commodities on it).
A {\em multicommodity transshipment} is a vector $x=(x^1,\dots,x^l)$
with $x^k\in\Z_+^t$ for all $k$ and $x^k_e$
the flow of commodity $k$ on edge $e$,
satisfying the capacity constraint $\sum_{k=1}^l x^k_e\leq u_e$
for each edge $e$ and demand constraint
$\sum_{e\in\delta^+(v)}x^k_e -\sum_{e\in\delta^-(v)}x^k_e=d^k_v$
for each vertex $v$ and commodity $k$
(with $\delta^+(v),\delta^-(v)$ the sets
of edges entering and leaving vertex $v$).

The cost of transshipment $x$ is defined as follows. There are cost functions
$f_e,g^k_e:\Z\rightarrow\Z$ for each edge and each edge-commodity pair.
The transshipment cost on edge $e$
is $f_e(\sum_{k=1}^l x^k_e)+\sum_{k=1}^l g^k_e(x^k_e)$ with the first
term being the value of $f_e$ on the combined flow of all
commodities on $e$ and the second term being the sum of costs that
depend on both the edge and the commodity. The total cost is
$$\sum_{e=1}^t\left(f_e\left(\sum_{k=1}^l x^k_e\right)
+\sum_{k=1}^l g^k_e(x^k_e)\right)\ .$$

Our results apply to cost functions which can
be standard linear or convex such as
$\alpha_e|\sum_{k=1}^l x^k_e|^{\beta_e}+
\sum_{k=1}^l \gamma^k_e |x^k_e|^{\delta^k_e}$
for some nonnegative integers $\alpha_e,\beta_e,\gamma^k_e,\delta^k_e$,
which take into account the increase in cost due to channel
congestion when subject to heavy traffic or communication load
(with the linear case obtained by $\beta_e=\delta^k_e$=1).

The theory of $n$-fold integer programming provides the first
polynomial time algorithms for the problem in two broad situations
discussed in \S \ref{tmctp} and \S \ref{tmtp}.

\subsubsection{The many-commodity transshipment problem}
\label{tmctp}

Here we consider the problem with {\em variable} number $l$ of commodities
over a fixed (but arbitrary) digraph - the so termed {\em many-commodity}
transshipment problem. This problem may seem at first very restricted:
however, even deciding if a feasible many-transshipment exists
(regardless of its cost) is NP-complete already over the complete
bipartite digraphs $K_{3,n}$ (oriented from one side to the other)
with only $3$ vertices on one side \cite{HOW2}; moreover, even
over the single tiny digraph $K_{3,3}$, the only solution
available to-date is the one given below via $n$-fold integer programming.

As usual, ${\hat f}$ and ${\hat g}$ denote the maximum
absolute values of the objective functions $f$ and $g$ over the feasible set.
It is usually easy to determine an upper bound on these values from the
problem data. For instance, in the special case of linear cost
functions $f$, $g$, bounds which are polynomial in the binary-encoding
length of the costs $\alpha_e$, $\gamma^k_e$, capacities $u$,
and demands $d^k_v$, are easily obtained by Cramer's rule.

We have the following theorem on (nonlinear) many-commodity transshipment.

\bt{Transshipment}{\bf \cite{HOW2}}
For every fixed digraph $G$ there is an algorithm that, given
$l$ commodity types, demand $d^k_v\in\Z$ for each commodity $k$
and vertex $v$, edge capacities $u_e\in\Z_+$, and convex costs
$f_e,g^k_e:\Z\rightarrow\Z$ presented by evaluation oracles,
solves in time polynomial in $l$ and $\l d^k_v,u_e,{\hat f},{\hat g}\r$,
the many-commodity transshipment problem,
\begin{eqnarray*}
&\min & \sum_e \left(f_e\left(\sum_{k=1}^l x^k_e\right)
+\sum_{k=1}^l g^k_e(x^k_e)\right) \\
&\mbox{\rm s.t.}& x^k_e\in\Z\,,\ \
\sum_{e\in\delta^+(v)}x^k_e -\sum_{e\in\delta^-(v)}x^k_e=d^k_v\,,\ \
\sum_{k=1}^l x^k_e\leq u_e\,,\ \ x^k_e\geq 0\ \ .
\end{eqnarray*}
\et
\bpr
Assume $G$ has $s$ vertices and $t$ edges and let $D$
be its $s\times t$ vertex-edge incidence matrix.
Let $f:\Z^t\rightarrow \Z$ and $g:\Z^{lt}\rightarrow \Z$ be the
separable convex functions defined by $f(y):=\sum_{e=1}^t f_e(y_e)$
with $y_e:=\sum_{k=1}^l x^k_e$
and $g(x):=\sum_{e=1}^t\sum_{k=1}^l g^k_e(x^k_e)$.
Let $x=(x^1,\dots,x^l)$ be the vector of variables
with $x^k\in\Z^t$ the flow of commodity $k$ for each $k$.
Then the problem can be rewritten in vector form as
\begin{equation*}
\min\left\{f\left(\sum_{k=1}^l x^k\right)+g\left(x\right)\ :\ x\in\Z^{lt}\,,\
Dx^k=d^k\,,\ \sum_{k=1}^l x^k\leq u\,,\ x\geq 0\right\}\ .
\end{equation*}
We can now proceed in two ways.

First way: extend the vector of variables to $x=(x^0,x^1,\dots,x^l)$
with $x^0\in\Z^t$ representing an additional slack commodity. Then the
capacity constraints become $\sum_{k=0}^l x^k=u$ and the cost function
becomes $f(u-x_0)+g(x^1,\dots,x^l)$ which is also separable convex.
Now let $A$ be the $(t,s)\times t$ bimatrix with first block
$A_1:=I_t$ the $t\times t$ identity matrix and second block $A_2:=D$.
Let $d^0:=Du-\sum_{k=1}^l d^k$ and let $b:=(u,d^0,d^1,\dots,d^l$).
Then the problem becomes the $(l+1)$-fold integer program
\begin{equation}\label{many-commodity_equation}
\min\left\{f\left(u-x^0\right)+g\left(x^1,\dots,x^l\right)
\ :\ x\in\Z^{(l+1)t}\,,\ A^{(l)}x=b\,,\ x\geq 0\right\}\ .
\end{equation}
By Theorem \ref{NFold2} this program
can be solved in polynomial time as claimed.

Second way: let $A$ be the $(0,s)\times t$ bimatrix with first block $A_1$ empty
and second block $A_2:=D$. Let $W$ be the $(t,0)\times t$ bimatrix with
first block $W_1:=I_t$ the $t\times t$ identity matrix and second block $W_2$ empty.
Let $b:=(d^1,\dots, d^l)$. Then the problem
is precisely the following $l$-fold integer program,
\begin{equation*}
\min\left\{f\left(W^{(l)}x\right)+g\left(x\right)
\ :\ x\in\Z^{lt}\,,\ A^{(l)}x=b\,,\ W^{(l)}x\leq u\,,\ x\geq 0\right\}\ .
\end{equation*}
By Theorem \ref{NFold5} this program
can be solved in polynomial time as claimed.
\epr

We also point out the following immediate corollary of Theorem \ref{Transshipment}.
\bc{CorollaryTransshipment}
For any fixed $s$, the (convex) many-commodity transshipment problem with variable
$l$ commodities on any $s$-vertex digraph is polynomial time solvable.
\ec

\subsubsection{The multicommodity transportation problem}
\label{tmtp}

Here we consider the problem with fixed (but arbitrary) number $l$
of commodities over any bipartite subdigraph of $K_{m,n}$
(oriented from one side to the other) - the so-called multicommodity
{\em transportation} problem - with fixed number $m$
of suppliers and {\em variable} number $n$ of consumers.
This is very natural in operations research
applications where few facilities serve many customers.
The problem is difficult even for $l=2$ commodities:
deciding if a feasible $2$-commodity transportation exists
(regardless of its cost) is NP-complete already over the complete
bipartite  digraphs
$K_{m,n}$ \cite{DO2}; moreover, even over the digraphs $K_{3,n}$
with only $m=3$ suppliers, the only available solution to-date
is the one given below via $n$-fold integer programming.

\begin{figure}[hbt]
\hskip-1.1cm
\includegraphics[scale=0.58]{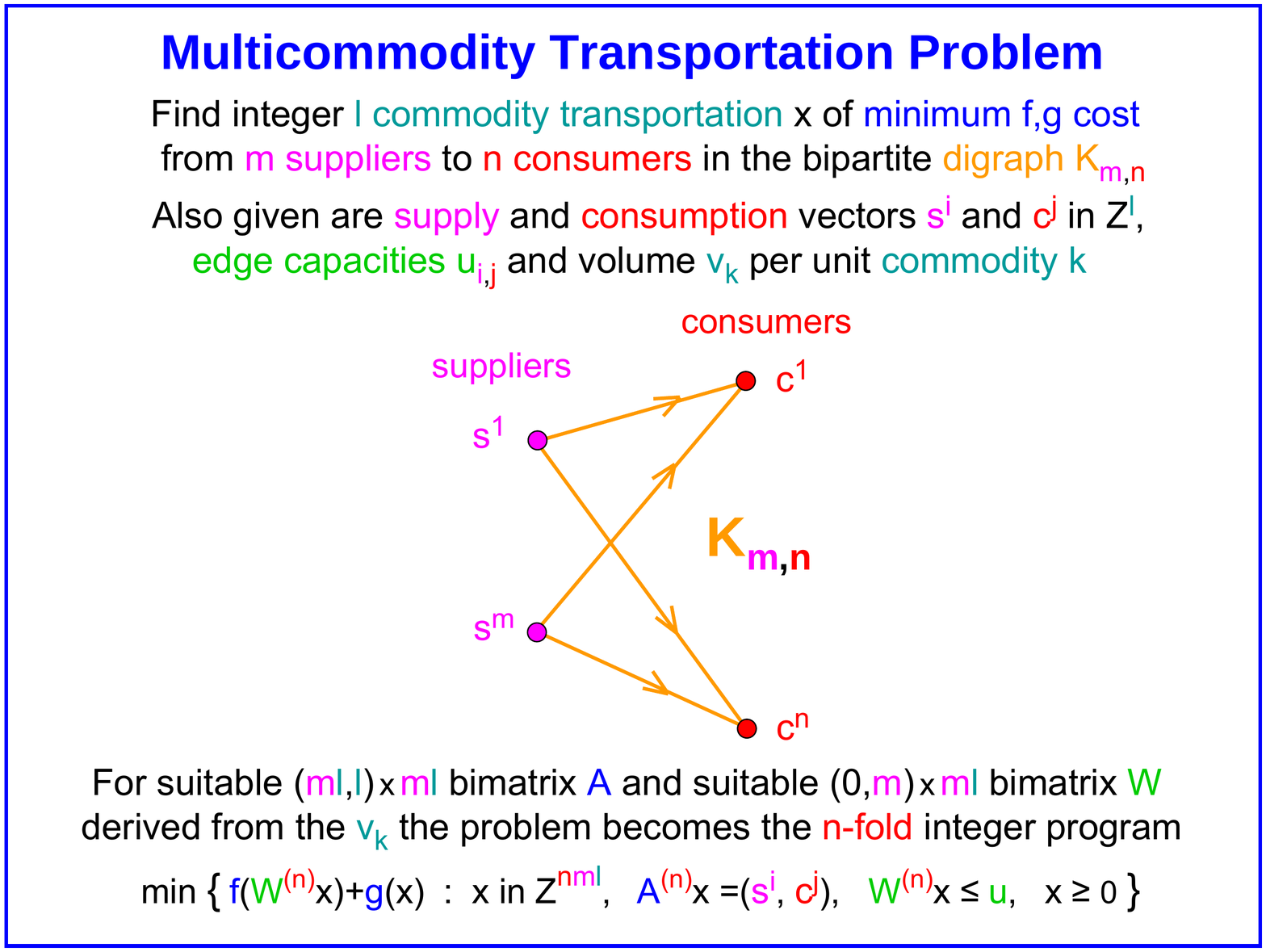}
\caption{Multicommodity Transportation Problem}
\label{multi-transportation-figure}
\end{figure}

This problem seems harder than the one discussed in the previous subsection
(with no seeming analog for non bipartite digraphs),
and the formulation below is more delicate. Therefore it is
convenient to change the labeling of the data a little bit as follows
(see Figure \ref{multi-transportation-figure}).
We now denote edges by pairs $(i,j)$ where $1\leq i\leq m$ is a supplier
and $1\leq j\leq n$ is a consumer. The demand vectors are now replaced
by (nonnegative) supply and consumption vectors: each supplier $i$
has a supply vector $s^i\in\Z_+^l$ with $s^i_k$ its supply in commodity $k$,
and each consumer $j$ has a consumption vector $c^j\in\Z_+^l$ with
$c^j_k$ its consumption in commodity $k$.
In addition, here each commodity $k$ may have its own volume
$v_k\in\Z_+$ per unit flow. A multicommodity transportation is now indexed as
$x=(x^1,\dots,x^n)$ with
$x^j=(x^j_{1,1},\dots,x^j_{1,l},\dots,x^j_{m,1},\dots,x^j_{m,l})$, where
$x^j_{i,k}$ is the flow of commodity $k$ from supplier $i$ to consumer $j$.
The capacity constraint on edge $(i,j)$ is
$\sum_{k=1}^l v_k x^j_{i,k}\leq u_{i,j}$ and the cost is
$f_{i,j}\left(\sum_{k=1}^l v_k x^j_{i,k}\right)+
\sum_{k=1}^l g^j_{i,k}\left(x^j_{i,k}\right)$
with $f_{i,j},g^j_{i,k}:\Z\rightarrow\Z$ convex.
As before, ${\hat f}$, ${\hat g}$ denote the maximum absolute
values of $f$, $g$ over the feasible set.

We assume below that the underlying digraph is $K_{m,n}$
(with edges oriented from suppliers to consumers), since the problem
over any subdigraph $G$ of $K_{m,n}$ reduces to that over $K_{m,n}$
by simply forcing $0$ capacity on all edges not present in $G$.

We have the following theorem on (nonlinear) multicommodity transportation.

\bt{Transportation}{\bf \cite{HOW2}}
For any fixed $l$ commodities, $m$ suppliers, and volumes
$v_k$, there is an algorithm that, given $n$, supplies and
demands $s^i,c^j\in\Z_+^l$, capacities $u_{i,j}\in\Z_+$, and convex
costs $f_{i,j},g^j_{i,k}:\Z\rightarrow\Z$ presented by evaluation oracles,
solves in time polynomial in $n$ and $\l s^i,c^j,u,{\hat f},{\hat g}\r$,
the multicommodity transportation problem,
\begin{eqnarray*}
&\min & \sum_{i,j} \left(f_{i,j}\left(\sum_k v_k x^j_{i,k}\right)
+\sum_{k=1}^l g^j_{i,k}\left(x^j_{i,k}\right)\right) \\
&\mbox{\rm s.t.}& x^j_{i,k}\in\Z\,,\ \
\sum_j x^j_{i,k}=s^i_k\,,\ \ \sum_i x^j_{i,k}=c^j_k\,,\ \
\sum_{k=1}^l v_k x^j_{i,k}\leq u_{i,j}\,,\ \ x^j_{i,k}\geq 0\ \ .
\end{eqnarray*}
\et
\bpr
Construct bimatrices $A$ and $W$ as follows.
Let $D$ be the $(l,0)\times l$ bimatrix with first block $D_1:=I_l$ and
second block $D_2$ empty. Let $V$ be the $(0,1)\times l$ bimatrix with
first block $V_1$ empty and second block $V_2:=(v_1,\dots,v_l)$.
Let $A$ be the $(ml,l)\times ml$ bimatrix with first block $A_1:=I_{ml}$
and second block $A_2:=D^{(m)}$. Let $W$ be the $(0,m)\times ml$ bimatrix
with first block $W_1$ empty and second block $W_2:=V^{(m)}$.
Let $b$ be the $(ml+nl)$-vector $b:=(s^1,\dots,s^m,c^1,\dots,c^n)$.

Let $f:\Z^{nm}\rightarrow \Z$ and $g:\Z^{nml}\rightarrow \Z$ be the
separable convex functions defined by $f(y):=\sum_{i,j}f_{i,j}(y_{i,j})$
with $y_{i,j}:=\sum_{k=1}^l v_k x^j_{i,k}$ and
$g(x):=\sum_{i,j}\sum_{k=1}^l g^j_{i,k}(x^j_{i,k})$.

Now note that $A^{(n)}x$ is an $(ml+nl)$-vector, whose first $ml$ entries
are the flows from each supplier of each commodity to all consumers,
and whose last $nl$ entries are the flows to each consumer of each
commodity from all suppliers. Therefore the supply and consumption
equations are encoded by $A^{(n)}x=b$. Next note that the
$nm$-vector $y=(y_{1,1},\dots,y_{m,1},\dots,y_{1,n},\dots,y_{m,n})$
satisfies $y=W^{(n)}x$. So the capacity constraints become $W^{(n)}x\leq u$
and the cost function becomes $f(W^{(n)}x)+g(x)$. Therefore, the
problem is precisely the following $n$-fold integer program,
\begin{equation*}
\min\left\{f\left(W^{(n)}x\right)+g\left(x\right)\ :\
x\in\Z^{nml}\,,\ A^{(n)}x=b\,,\ W^{(n)}x\leq u\,,\ x\geq 0\right\}\ .
\end{equation*}
By Theorem \ref{NFold5} this program can be solved in polynomial time as claimed.
\epr

\section{Theory}
\label{t}
In \S \ref{gbanip} we define Graver bases of integer matrices and show that
they can be used to solve linear and nonlinear integer programs in polynomial time.
In \S \ref{nfip} we show that Graver bases of $n$-fold products can be computed
in polynomial time and, incorporating the results of \S \ref{gbanip},
prove our Theorems \ref{NFold1}--\ref{NFold5} that establish the polynomial
time solvability of linear and nonlinear $n$-fold integer programming.

To simplify the presentation, and since the feasible set in most
applications is finite or can be made finite by more careful
modeling, whenever an algorithm detects that the feasible set
is infinite, it simply stops. So, throughout our discussion, an algorithm
is said to {\em solve a (nonlinear) integer programming problem} if it either
finds an optimal solution $x$ or concludes that the feasible set is infinite or empty.

As noted in the introduction, any nonlinear function $f$ involved
is presented either by a mere comparison oracle that
for any two vectors $x,y$ can answer whether or not $f(x)\leq f(y)$, or by an
evaluation oracle that for any vector $x$ can return $f(x)$.

\subsection{Graver Bases and Nonlinear Integer Programming}
\label{gbanip}

The {\em Graver basis} is a fundamental object in the theory of integer programming
which was introduced by J. Graver already back in 1975 \cite{Gra}.
However, only very recently, in the series of papers \cite{DHOW,DHORW,HOW1},
it was established that the Graver basis can be used to solve linear
(as well as nonlinear) integer programming problems in polynomial time.
In this subsection we describe these important new developments.

\subsubsection{Graver bases}
\label{gb}

We begin with the definition of the Graver basis and some of its basic properties.
Throughout this subsection let $A$ be an integer $m\times n$ matrix.
The \emph{lattice}\index{lattice} of $A$ is the set
$\L(A):=\{x\in\Z^n\,:\, Ax=0\}$ of integer vectors in its kernel.
We use $\L^*(A):=\{x\in\Z^n\,:\, Ax=0,\ x\neq 0\}$ to denote
the set of nonzero elements in $\L(A)$. We use a partial order $\sqsubseteq$
on $\R^n$ which extends the usual coordinate-wise partial order $\leq$ on the
nonnegative orthant $\R^n_+$ and is defined as follows. For two vectors
$x,y\in\R^n$ we write $x\sqsubseteq y$ and say that $x$ is {\em conformal} to $y$
if $x_iy_i\geq 0$ and $|x_i|\leq |y_i|$ for $i=1,\ldots,n$, that is, $x$ and $y$
lie in the same orthant of $\R^n$ and each component of $x$ is bounded by
the corresponding component of $y$ in absolute value. A suitable extension
of the classical lemma of Gordan \cite{Gor} implies that every subset of $\Z^n$
has finitely many $\sqsubseteq$-minimal elements. We have the following
fundamental definition.

\bd{GraverBasisDefinition}{\bf \cite{Gra}}
The {\em Graver basis} of an integer matrix $A$ is defined to be the
finite set $\G(A)\subset\Z^n$ of $\sqsubseteq$-minimal elements
in $\L^*(A)=\{x\in\Z^n\,:\, Ax=0,\ x\neq 0\}$.
\ed
Note that $\G(A)$ is centrally symmetric, that is, $g\in\G(A)$ if and only
if $-g\in\G(A)$. For instance, the Graver basis of the $1\times 3$
matrix $A:=[1\,\ 2\,\ 1]$ consists of $8$ elements,
$$\G(A)\ =\ \pm\left\{(2,-1,0),(0,-1,2),(1,0,-1),(1,-1,1)\right\}\ .$$
Note also that the Graver basis may contain elements, such as $(1,-1,1)$
in the above small example, whose support involves linearly dependent columns
of $A$. So the cardinality of the Graver basis cannot be bounded in
terms of $m$ and $n$ only and depends on the entries of $A$ as well.
Indeed, the Graver basis is typically exponential and cannot be written
down, let alone computed, in polynomial time. But, as we will show in the
next section, for $n$-fold products it can be computed efficiently.

\vskip.2cm
A finite sum $u:=\sum_i v_i$ of vectors in $\R^n$ is called {\em conformal}
if all summands lie in the same orthant and hence $v_i\sqsubseteq u$ for
all $i$. We start with a simple lemma.
\bl{ConformalGraverSumWeakForm}
Any $x\in\L^*(A)$ is a conformal sum $x=\sum_i g_i$ of Graver basis elements
$g_i\in\G(A)$, with some elements possibly appearing more than once in the sum.
\el
\bpr
We use induction on the well partial order $\sqsubseteq$.
Consider any $x\in\L^*(A)$. If it is $\sqsubseteq$-minimal in $\L^*(A)$
then $x\in\G(A)$ by definition of the Graver basis and we are done.
Otherwise, there is an element $g\in\G(A)$ such that $g\sqsubset x$.
Set $y:=x-g$. Then $y\in\L^*(A)$ and $y\sqsubset x$, so by induction
there is a conformal sum $y=\sum_i g_i$ with $g_i\in\G(A)$ for all $i$.
Now $x=g+\sum_i g_i$ is a conformal sum of $x$.
\epr
We now provide a stronger form of Lemma \ref{ConformalGraverSumWeakForm}
which basically follows from the integer analogs of Carath\'eodory's theorem
established in \cite{CFS} and \cite{Seb}.
\bl{ConformalGraverSum}
Any $x\in\L^*(A)$ is a conformal sum $x=\sum_{i=1}^t\lambda_i g_i$ involving
$t\leq 2n-2$ Graver basis elements $g_i\in\G(A)$
with nonnegative integer coefficients $\lambda_i\in\Z_+$.
\el
\bpr
We prove the slightly weaker bound $t\leq 2n-1$ from \cite{CFS}.
A proof of the stronger bound can be found in \cite{Seb}.
Consider any $x\in\L^*(A)$ and let $g_1,\dots,g_s$ be all elements
of $\G(A)$ lying in the same orthant as $x$. Consider the linear program
\begin{equation}\label{ConformalLP}
\max\left\{\sum_{i=1}^s\lambda_i\ :\
x=\sum_{i=1}^s\lambda_i g_i\,,\ \ \lambda_i\in\R_+\right\}\ .
\end{equation}
By Lemma \ref{ConformalGraverSumWeakForm} the point $x$ is a nonnegative linear
combination of the $g_i$ and hence the program (\ref{ConformalLP}) is feasible.
Since all $g_i$ are nonzero and in the same orthant as $x$,
program (\ref{ConformalLP}) is also bounded. As is well known, it then
has a {\em basic} optimal solution, that is, an optimal solution
$\lambda_1,\dots,\lambda_s$ with at most $n$ of the $\lambda_i$ nonzero. Let
$$y\ :=\ \sum(\lambda_i-\lfloor\lambda_i\rfloor)g_i
\ =\ x-\sum\lfloor\lambda_i\rfloor g_i\ .$$
If $y=0$ then $x=\sum\lfloor\lambda_i\rfloor g_i$ is a conformal sum of at
most $n$ of the $g_i$ and we are done. Otherwise, $y\in\L^*(A)$ and $y$ lies
in the same orthant as $x$, and hence, by Lemma \ref{ConformalGraverSumWeakForm}
again, $y=\sum_{i=1}^s\mu_i g_i$ with all $\mu_i\in\Z_+$.
Then $x=\sum(\mu_i+\lfloor\lambda_i\rfloor)g_i$ and hence,
since the $\lambda_i$ form an optimal solution to (\ref{ConformalLP}), we have
$\sum(\mu_i+\lfloor\lambda_i\rfloor)\leq\sum\lambda_i$. Therefore
$\sum \mu_i\leq\sum(\lambda_i-\lfloor\lambda_i\rfloor)<n$ with the last inequality
holding since at most $n$ of the $\lambda_i$ are nonzero. Since the $\mu_i$ are integer,
at most $n-1$ of them are nonzero. So $x=\sum(\mu_i+\lfloor\lambda_i\rfloor)g_i$
is a conformal sum of $x$ involving at most $2n-1$ of the $g_i$.
\epr

\vskip.2cm
The Graver basis also enables to check the finiteness of a feasible integer program.
\bl{GraverFiniteness}
Let $\G(A)$ be the Graver basis of matrix $A$ and let $l,u\in\Z_{\infty}^n$.
If there is some $g\in\G(A)$ satisfying $g_i\leq 0$ whenever $u_i<\infty$
and $g_i\geq 0$ whenever $l_i>-\infty$ then every set of the form
$S:=\{x\in\Z^n\,:\,Ax=b\,,\ l\leq x\leq u\}$ is either empty or infinite,
whereas if there is no such $g$, then every set $S$ of this form is finite.
Clearly, the existence of such $g$ can be checked in time polynomial in
$\l \G(A),l,u \r$.
\el
\bpr
First suppose there exists such $g$. Consider any such $S$.
Suppose $S$ contains some point $x$. Then for all $\lambda\in\Z_+$ we
have $l\leq x+\lambda g\leq u$ and $A(x+\lambda g)=Ax=b$ and hence
$x+\lambda g\in S$, so $S$ is infinite. Next suppose $S$ is infinite.
Then the polyhedron $P:=\{x\in\R^n\,:\,Ax=b\,,\ l\leq x\leq u\}$
is unbounded and hence, as is well known, has a recession vector,
that is, a nonzero $h$, which we may assume to be integer,
such that $x+\alpha h\in P$ for all $x\in P$ and $\alpha\geq 0$.
This implies that $h\in\L^*(A)$ and that $h_i\leq 0$ whenever $u_i<\infty$
and $h_i\geq 0$ whenever $l_i>-\infty$. So $h$ is a conformal sum
$h=\sum g_i$ of vectors $g_i\in\G(A)$, each of which also
satisfies $g_i\leq 0$ whenever $u_i<\infty$ and $g_i\geq 0$
whenever $l_i>-\infty$, providing such $g$.
\epr

\subsubsection{Separable convex integer minimization}
\label{scim}

In this subsection we consider the following nonlinear integer minimization problem
\begin{equation}\label{SeparableProblem}
\min\{f(x)\ :\ x\in\Z^n,\ Ax=b,\ l\leq x\leq u\}\ ,
\end{equation}
where $A$ is an integer $m\times n$ matrix, $b\in\Z^m$, $l,u\in\Z_{\infty}^n$,
and $f:\Z^n\rightarrow\Z$ is a separable convex function, that is,
$f(x)=\sum_{j=1}^n f_j(x_j)$ with $f_j:\Z\rightarrow\Z$
a univariate convex function for all $j$. We prove a sequence of
lemmas and then combine them to show that the Graver basis of $A$
enables to solve this problem in polynomial time.

\vskip.2cm
We start with two simple lemmas about univariate convex functions. The first
lemma establishes a certain {\em supermodularity} property of such functions.
\bl{Supermodular}
Let $f:\R\rightarrow\R$ be a univariate convex function, let $r$ be a real number,
and let $s_1,\dots,s_m$ be real numbers satisfying $s_is_j\geq 0$ for all $i,j$.
Then we have
$$f\left(r+\sum_{i=1}^m s_i\right)-f(r)\ \geq\
\sum_{i=1}^m \left(f(r+s_i)-f(r)\right)\ .$$
\el
\bpr
We use induction on $m$. The claim holding trivially for $m=1$, consider $m>1$.
Since all nonzero $s_i$ have the same sign, $s_m=\lambda\sum_{i=1}^m s_i$
for some $0\leq\lambda\leq1$. Then
$$r+s_m=(1-\lambda) r + \lambda\left(r+\sum_{i=1}^m s_i\right)\,,\quad
r+\sum_{i=1}^{m-1}s_i=\lambda r + (1-\lambda)\left(r+\sum_{i=1}^m s_i\right)\ ,$$
and so the convexity of $f$ implies
\begin{eqnarray*}
f(r+s_m) & + & f\left(r+\sum_{i=1}^{m-1}s_i\right) \\
&\leq & (1-\lambda)f(r) + \lambda f\left(r+\sum_{i=1}^m s_i\right)
 + \lambda f(r) + (1-\lambda)f\left(r+\sum_{i=1}^m s_i\right) \\
&\ & \hskip7.1cm =\ f(r)\ +\ f\left(r+\sum_{i=1}^m s_i\right).
\end{eqnarray*}
Subtracting $2f(r)$ from both sides and applying induction, we obtain, as claimed,
\begin{eqnarray*}
f\left(r+\sum_{i=1}^m s_i\right) & - & f(r) \\
&\geq & f(r+s_m)-f(r) + f\left(r+\sum_{i=1}^{m-1}s_i\right)-f(r) \\
& \ & \hskip3.15cm\geq\ \sum_{i=1}^m \left(f(r+s_i)-f(r)\right)\ .
\mepr
\end{eqnarray*}

The second lemma shows that univariate convex functions can be minimized
efficiently over an interval of integers using repeated bisections.
\bl{UnivariateGreedyAugmentation}
There is an algorithm that, given any two integer numbers $r\leq s$ and any univariate
convex function $f:\Z\rightarrow\R$ given by a comparison oracle, solves in time
polynomial in $\l r,s \r$ the following univariate integer minimization problem,
$$\min\left\{\,f(\lambda)\ :\ \lambda\in\Z\,,\ \ r\leq\lambda\leq s\,\right\} .$$
\el
\bpr
If $r=s$ then $\lambda:=r$ is optimal. Assume then $r\leq s-1$.
Consider the integers
$$r\ \leq\ \left\lfloor {r+s\over2}\right\rfloor\ <\
\left\lfloor {r+s\over2}\right\rfloor+1\ \leq\  s\ .$$
Use the oracle of $f$ to compare
$f\left(\left\lfloor {r+s\over2}\right\rfloor\right)$
and $f\left(\left\lfloor {r+s\over2}\right\rfloor+1\right)$.
By the convexity of $f$:
$$\begin{array}{ll}
    f\left(\lfloor {r+s\over2}\rfloor\right)=
    f\left(\lfloor {r+s\over2}\rfloor+1\right)\ \Rightarrow\
      & \hbox{$\lambda:=\lfloor {r+s\over2}\rfloor$ is a minimum of $f$;} \\
    f\left(\lfloor {r+s\over2}\rfloor\right)<
    f\left(\lfloor {r+s\over2}\rfloor+1\right)\ \Rightarrow\
      & \hbox{the minimum of $f$ is in the interval
              $[r,\left\lfloor {r+s\over2}\right\rfloor]$;} \\
    f\left(\lfloor {r+s\over2}\rfloor\right)>
    f\left(\lfloor {r+s\over2}\rfloor+1\right)\ \Rightarrow\
      & \hbox{the minimum of $f$ is in the interval
              $[\left\lfloor {r+s\over2}\right\rfloor+1,s]$.} \\
\end{array}$$
Thus, we either obtain the optimal point, or bisect the interval
$[r,s]$ and repeat. So in $O(\log (s-r))=O(\l r,s\r)$ bisections
we find an optimal solution $\lambda\in\Z\cap[r,s]$.
\epr

\begin{figure}[hbt]
\hskip-1.1cm
\includegraphics[scale=0.58]{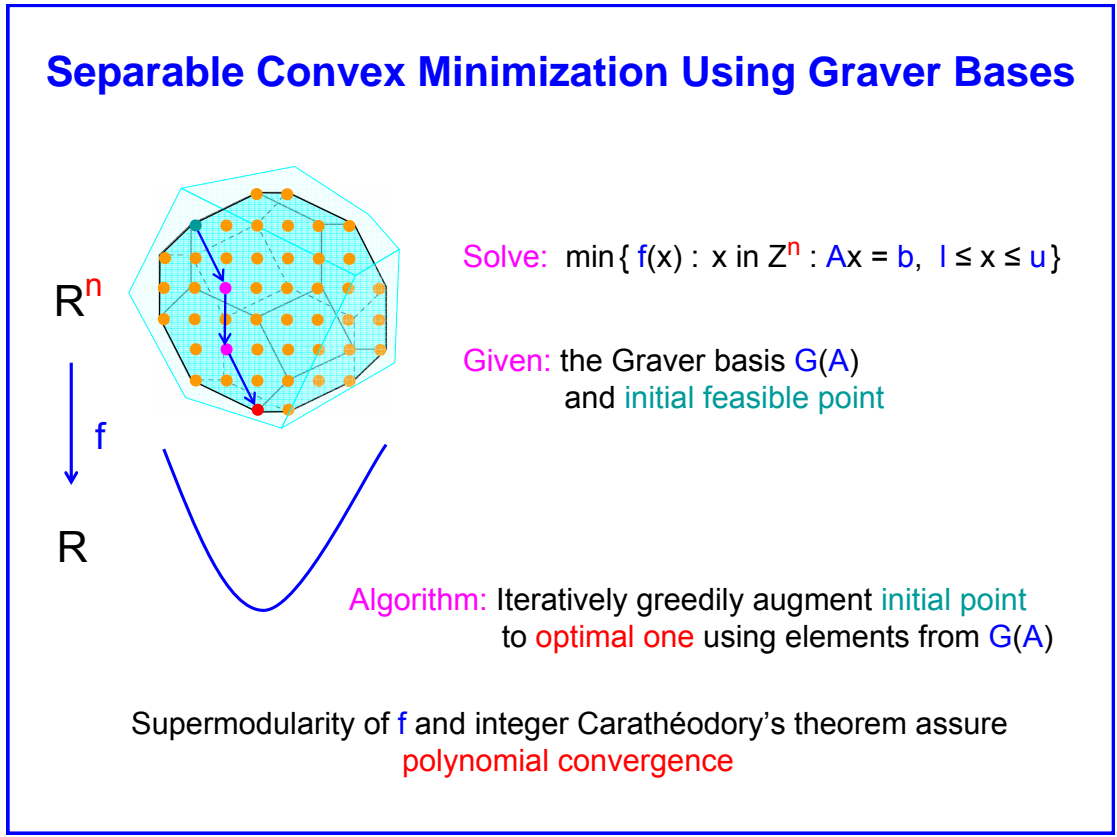}
\caption{Separable Convex Minimization Using Graver Bases}
\label{convex_minimization_figure}
\end{figure}

The next two lemmas extend Lemmas \ref{Supermodular}
and \ref{UnivariateGreedyAugmentation}.
The first lemma shows the supermodularity of
separable convex functions with respect to conformal sums.
\bl{SeparableConvexConformal}
Let $f:\R^n\rightarrow\R$ be any separable convex function, let $x\in\R^n$ be any point,
and let $\sum g_i$ be any conformal sum in $\R^n$. Then the following inequality holds,
$$f\left(x+\sum g_i\right)-f(x)\ \geq\ \sum \left(f\left(x+g_i\right)-f(x)\right)\ .$$
\el
\bpr
Let $f_j:\R\rightarrow\R$ be univariate convex functions such that
$f(x)=\sum_{j=1}^n f_j(x_j)$. Consider any $1\leq j\leq n$.
Since $\sum g_i$ is a conformal sum, we have $g_{i,j}g_{k,j}\geq 0$
for all $i,k$ and so, setting $r:=x_j$ and $s_i:=g_{i,j}$ for all $i$,
Lemma \ref{Supermodular} applied to $f_j$ implies
\begin{equation}\label{separable_equation}
f_j\left(x_j+\sum_i g_{i,j}\right)-f_j(x_j)\ \geq\
\sum_i \left(f_j\left(x_j+g_{i,j}\right)-f_j(x_j)\right)\ .
\end{equation}
Summing the equations (\ref{separable_equation}) for $j=1,\dots,n$,
we obtain the claimed inequality.
\epr
The second lemma concerns finding a best improvement step in a given direction.
\bl{GreedyAugmentation}
There is an algorithm that, given bounds $l,u\in\Z_{\infty}^n$,
direction $g\in\Z^n$, point $x\in\Z^n$ with $l\leq x\leq u$,
and convex function $f:\Z^n\rightarrow\R$ presented by comparison oracle,
solves in time polynomial in $\l l,u,g,x \r$, the univariate problem,
\begin{equation}\label{GreedyAugmentationEquation}
\min\{f(x+\lambda g)\ :\ \lambda\in\Z_+\,,\ l\leq x+\lambda g\leq u\}
\end{equation}
\el
\bpr
Let $S:=\{\lambda\in\Z_+\,:\,l\leq x+\lambda g\leq u\}$ be the feasible set
and let $s:=\sup S$, which is easy to determine. If $s=\infty$ then conclude
that $S$ is infinite and stop. Otherwise, $S=\{0,1,\dots,s\}$ and the problem
can be solved by the algorithm of Lemma \ref{UnivariateGreedyAugmentation}
minimizing the univariate convex function $h(\lambda):=h(x+\lambda g)$ over $S$.
\epr
We can now show that the Graver basis of $A$ allows to solve problem
(\ref{SeparableProblem}) in polynomial time, provided we are given an initial feasible
point to start with. We will later show how to find such an initial point as well.
As noted in the introduction, ${\hat f}$ below denotes the maximum value of $|f(x)|$
over the feasible set (which need not be part of the input). An outline
of the algorithm is provided in Figure \ref{convex_minimization_figure}.
\bl{SeparableConvexMinimizationLemma}
There is an algorithm that, given an integer $m\times n$ matrix $A$,
its Graver basis $\G(A)$, vectors $l,u\in\Z_{\infty}^n$ and $x\in\Z^n$ with
$l\leq x\leq u$, and separable convex function $f:\Z^n\rightarrow\Z$ presented by a
comparison oracle, solves the integer program
\begin{equation}\label{WeakSeparableConvexEquation}
\min\{f(z)\ :\ z\in\Z^n\,,\ Az=b\,,\ l\leq z\leq u\}\ ,\quad b:=Ax\ ,
\end{equation}
in time polynomial in the binary-encoding length $\l\G(A),l,u,x,{\hat f}\r$
of the data.
\el
\bpr
First, apply the algorithm of Lemma \ref{GraverFiniteness} to $\G(A)$
and $l,u$ and either detect that the feasible set is infinite and stop,
or conclude it is finite and continue. Next produce a sequence of
feasible points $x_0,x_1,\ldots,x_s$ with $x_0:=x$ the given input point,
as follows. Having obtained $x_k$, solve the univariate minimization problem
\begin{equation}\label{GreedyEquation}
\min\{f(x_k+\lambda g)\ :\ \lambda\in\Z_+\,,\
g\in \G(A)\,,\ l\leq x_k+\lambda g\leq u\,\}
\end{equation}
by applying the algorithm of Lemma \ref{GreedyAugmentation} for each $g\in \G(A)$.
If the minimal value in (\ref{GreedyEquation})
satisfies $f(x_k+\lambda g)<f(x_k)$ then set $x_{k+1}:=x_k+\lambda g$
and repeat, else stop and output the last point $x_s$ in the sequence.
Now, $Ax_{k+1}=A(x_k+\lambda g)=Ax_k=b$ by induction on $k$, so
each $x_k$ is feasible.
Since the feasible set is finite and the $x_k$ have decreasing
objective values and hence distinct, the algorithm terminates.

We now show that the point $x_s$ output by the algorithm is optimal.
Let $x^*$ be any optimal solution to (\ref{WeakSeparableConvexEquation}).
Consider any point $x_k$ in the sequence and suppose it is not optimal.
We claim that a new point $x_{k+1}$ will be produced and will satisfy
\begin{equation}\label{FastConvergence}
f(x_{k+1})-f\left(x^*\right)\ \leq\ {2n-3\over 2n-2}\left(f(x_k)-f(x^*)\right)
\end{equation}
By Lemma \ref{ConformalGraverSum}, we can write the difference
$x^*-x_k=\sum_{i=1}^t\lambda_i g_i$ as conformal sum involving $1\leq t\leq 2n-2$ elements
$g_i\in \G(A)$ with all $\lambda_i\in\Z_+$. By Lemma \ref{SeparableConvexConformal},
$$f(x^*)-f\left(x_k\right)\ =\ f\left(x_k+\sum_{i=1}^t\lambda_i g_i\right)-f(x_k)
\ \geq\ \sum_{i=1}^t\left(f\left(x_k+\lambda_i g_i\right)-f(x_k)\right)\ .$$
Adding $t\left(f(x_k)-f(x^*)\right)$ on both sides and rearranging terms we obtain
$$\sum_{i=1}^t\left(f\left(x_k+\lambda_i g_i\right)-f(x^*)\right)
\ \leq\ (t-1)\left(f(x_k)-f(x^*)\right)\ .$$
Therefore there is some summand on the left-hand side satisfying
$$f\left(x_k+\lambda_i g_i\right)-f(x^*)\ \leq\ {t-1\over t}\left(f(x_k)-f(x^*)\right)
\ \leq\ {2n-3\over 2n-2}\left(f(x_k)-f(x^*)\right)\ .$$
So the point $x_k+\lambda g$ attaining minimum in (\ref{GreedyEquation}) satisfies
$$f(x_k+\lambda g)-f(x^*)\ \leq\ f\left(x_k+\lambda_i g_i\right)-f(x^*)
\ \leq\ {2n-3\over 2n-2}\left(f(x_k)-f(x^*)\right)$$ and so indeed
$x_{k+1}:=x_k+\lambda g$ will be produced and will satisfy (\ref{FastConvergence}).
This shows that the last point $x_s$ produced and output by the algorithm is indeed optimal.

We proceed to bound the number $s$ of points. Consider any $i<s$ and the
intermediate non optimal point $x_i$ in the sequence produced by the algorithm.
Then $f(x_i)>f(x^*)$ with both values integer, and so
repeated use of (\ref{FastConvergence}) gives
\begin{eqnarray*}
1\leq f(x_i)-f(x^*) & = & \prod_{k=0}^{i-1}
{{f(x_{k+1})-f(x^*)}\over{f(x_k)-f(x^*)}}\left(f(x)-f(x^*)\right) \\
& \leq & \left({2n-3\over 2n-2}\right)^i\left(f(x)-f(x^*)\right)
\end{eqnarray*}
and therefore
$$i\ \leq\ \left(\log{2n-2\over 2n-3}\right)^{-1}\log\left(f(x)-f(x^*)\right)\ .$$
Therefore the number $s$ of points produced by the algorithm is at most one unit
larger than this bound, and using a simple bound on the logarithm, we obtain
$$s\ =\ O\left(n \log(f(x)-f(x^*))\right)\ .$$
Thus, the number of points produced and the total running time are polynomial.
\epr
Next we show that Lemma \ref{SeparableConvexMinimizationLemma}
can also be used to find an initial feasible point for the given
integer program or assert that none exists in polynomial time.
\bl{Feasibility}
There is an algorithm that, given integer $m\times n$ matrix $A$,
its Graver basis $\G(A)$, $l,u\in\Z_{\infty}^n$, and $b\in\Z^m$,
either finds an $x\in\Z^n$ satisfying $l\leq x\leq u$ and $Ax=b$
or asserts that none exists, in time which is polynomial in $\l A,\G(A),l,u,b\r$.
\el
\bpr
Assume that $l\leq u$ and that $l_i<\infty$ and
$u_j>-\infty$ for all $j$, since otherwise there is no feasible point.
Also assume that there is no $g\in\G(A)$ satisfying $g_i\leq 0$
whenever $u_i<\infty$ and $g_i\geq 0$ whenever $l_i>-\infty$, since
otherwise $S$ is empty or infinite by Lemma \ref{GraverFiniteness}.
Now, either detect there is no integer solution
to the system of equations $Ax=b$ (without the lower
and upper bound constraints) and stop, or determine some such solution
${\hat x}\in\Z^n$ and continue; it is well known that this can be done in
polynomial time, say, using the Hermite normal form of $A$,
see \cite{Sch}. Next define a separable convex function on $\Z^n$ by
$f(x):=\sum_{j=1}^n f_j(x_j)$ with
$$f_j(x_j)\ :=\
\left\{
  \begin{array}{ll}
    l_j-x_j, & \hbox{if $x_j< l_j$} \\
    0, & \hbox{if $l_j\leq x_j\leq u_j$} \\
    x_j-u_j, & \hbox{if $x_j> u_j$}
  \end{array}
\right.,\quad\quad j=1,\dots,n
$$
and extended lower and upper bounds
$${\hat l}_j\ :=\ \min\{l_j,{\hat x}_j\}\,,\quad\quad
{\hat u}_j\ :=\ \max\{u_j,{\hat x}_j\}\,,\quad\quad j=1,\dots,n\ .$$
Consider the auxiliary separable convex integer program
\begin{equation}\label{auxiliary_program}
\min\{f(z)\ :\ z\in\Z^n\,,\ Az=b\,,\ {\hat l}\leq z\leq {\hat u}\}
\end{equation}
First note that ${\hat l}_j>-\infty$ if and only if $l_j>-\infty$
and ${\hat u}_j<\infty$ if and only if $u_j<\infty$. Therefore
there is no $g\in\G(A)$ satisfying $g_i\leq 0$ whenever
${\hat u}_i<\infty$ and $g_i\geq 0$ whenever ${\hat l}_i>-\infty$
and hence the feasible set of (\ref{auxiliary_program}) is finite
by Lemma \ref{GraverFiniteness}. Next note that ${\hat x}$ is feasible
in (\ref{auxiliary_program}). Now apply the algorithm of Lemma
\ref{SeparableConvexMinimizationLemma} to (\ref{auxiliary_program})
and obtain an optimal solution $x$. Note that this can be done in
polynomial time since the binary length of ${\hat x}$ and therefore
also of $\hat l$, $\hat u$ and of the maximum value $\hat f$
of $|f(x)|$ over the feasible set of (\ref{auxiliary_program})
are polynomial in the length of the data.

Now note that every point $z\in S$ is feasible in (\ref{auxiliary_program}),
and every point $z$ feasible in (\ref{auxiliary_program})
satisfies $f(z)\geq 0$ with equality if and only if $z\in S$.
So, if $f(x)>0$ then the original set $S$ is empty,
whereas if $f(x)=0$ then $x\in S$ is a feasible point.
\epr
We are finally in position, using Lemmas \ref{SeparableConvexMinimizationLemma}
and \ref{Feasibility}, to show that the Graver basis allows to solve the
nonlinear integer program (\ref{SeparableProblem}) in polynomial time.
As usual, ${\hat f}$ is the maximum of $|f(x)|$ over
the feasible set and need not be part of the input.
\bt{Graver2}{\bf \cite{HOW1}}
There is an algorithm that, given integer $m\times n$ matrix $A$,
its Graver basis $\G(A)$, $l,u\in\Z_{\infty}^n$, $b\in\Z^m$,
and separable convex $f:\Z^n\rightarrow\Z$ presented
by comparison oracle, solves in time polynomial in
$\l A,\G(A),l,u,b,{\hat f}\r$ the problem
$$\min\{f(x)\ :\ x\in\Z^n\,,\ Ax=b\,,\ l\leq x\leq u\}\ .$$
\et
\bpr
First, apply the polynomial time algorithm of Lemma \ref{Feasibility} and either
conclude that the feasible set is infinite or empty and stop,
or obtain an initial feasible point and continue. Next, apply
the polynomial time algorithm of Lemma \ref{SeparableConvexMinimizationLemma}
and either conclude that the feasible set is infinite
or obtain an optimal solution.
\epr

\subsubsection{Specializations and extensions}
\label{lascm}

\vskip.5cm
\subsubsection*{\sc Linear integer programming}
\label{lip}

Any linear function $wx=\sum_{i=1}^n w_ix_i$ is separable convex.
Moreover, an upper bound on $|wx|$ over the feasible set (when finite),
which is polynomial in the binary-encoding length of the data,
readily follows from Cramer's rule. Therefore we obtain, as an
immediate special case of Theorem \ref{Graver2},
the following important result, asserting that Graver
bases enable the polynomial time solution of linear integer programming.
\bt{Graver1}{\bf \cite{DHOW}}
There is an algorithm that, given an integer $m\times n$ matrix $A$,
its Graver basis $\G(A)$, $l,u\in\Z_{\infty}^n$, $b\in\Z^m$, and $w\in\Z^n$,
solves in time which is polynomial in $\l A,\G(A),l,u,b,w\r$,
the following linear integer programming problem,
$$\min\{wx\ :\ x\in\Z^n\,,\ Ax=b\,,\ l\leq x\leq u\}\ .$$
\et

\vskip.5cm
\subsubsection*{\sc Distance minimization}
\label{dm}

Another useful special case of Theorem \ref{Graver2} which is natural in
various applications such as image processing, tomography, communication,
and error correcting codes, is the following result, which asserts
that the Graver basis enables to determine a feasible point
which is $l_p$-closest to a given desired goal point in polynomial time.
\bt{Graver3}{\bf \cite{HOW1}}
There is an algorithm that, given integer $m\times n$ matrix $A$, its Graver
basis $\G(A)$, positive integer $p$, vectors $l,u\in\Z_{\infty}^n$,
$b\in\Z^m$, and ${\hat x}\in\Z^n$, solves in time polynomial in $p$
and $\l A,\G(A),l,u,b,{\hat x}\r$, the distance minimization problem
\begin{equation}\label{distance_equation_second}
\min\,\{\|x-{\hat x}\|_p\ :\ x\in\Z^n,\ Ax=b,\ l\leq x\leq u\}\ .
\end{equation}
For $p=\infty$ the problem (\ref{distance_equation_second})
can be solved in time polynomial in $\l A,\G(A),l,u,b,{\hat x}\r$.
\et
\bpr
For finite $p$ apply the algorithm of Theorem \ref{Graver2}
taking $f$ to be the $p$-th power $\|x-{\hat x}\|_p^p$ of
the $l_p$ distance. If the feasible set is nonempty and finite
(else the algorithm stops) then the maximum value $\hat f$ of $|f(x)|$
over it is polynomial in $p$ and $\l A,l,u,b,{\hat x}\r$,
and hence an optimal solution can be found in polynomial time.

Consider $p=\infty$. Using Cramer's rule it is easy to compute
an integer $\rho$ with $\l\rho\r$ polynomially bounded in
$\l A,l,u,b\r$ that, if the feasible set is finite,
provides an upper bound on $\|x\|_\infty$ for any feasible $x$ .
Let $q$ be a positive integer satisfying
$$q\ >\ {{\log n}\over{\log(1+(2\rho)^{-1})}}\ .$$
Now apply the algorithm of the first paragraph above for the $l_q$ distance.
Assuming the feasible set is nonempty and finite (else the algorithm stops)
let $x^*$ be the feasible point which minimizes the $l_q$ distance
to $\hat x$ obtained by the algorithm. We claim that it also minimizes
the $l_\infty$ distance to $\hat x$ and hence is the desired optimal solution.
Consider any feasible point $x$. By standard inequalities between the
$l_\infty$ and $l_q$ norms,
$$\|x^*-{\hat x}\|_\infty\ \leq\ \|x^*-{\hat x}\|_q \ \leq\
\|x-{\hat x}\|_q\ \leq\ n^{1\over q} \|x-{\hat x}\|_\infty\ .$$
Therefore
$$\|x^*-{\hat x}\|_\infty-\|x-{\hat x}\|_\infty\ \leq\
(n^{1\over q}-1)\|x-{\hat x}\|_\infty\ \leq\ (n^{1\over q}-1)2\rho\ <\ 1\ ,$$
where the last inequality holds by the choice of $q$.
Since $\|x^*-{\hat x}\|_\infty$ and $\|x-{\hat x}\|_\infty$
are integers we find that
$\|x^*-{\hat x}\|_\infty\leq \|x-{\hat x}\|_\infty$.
This establishes the claim.
\epr
In particular, for all positive $p\in\Z_\infty$, using the Graver basis we can solve
$$\min\,\{\|x\|_p\ :\ x\in\Z^n,\ A x=b,\ l\leq x\leq u\}\ ,$$
which for $p=\infty$ is equivalent to the min-max integer program
$$\min\left\{\max\{|x_i|\,:\, i=1,\dots,n\}
\ :\ x\in\Z^n,\ Ax=b,\ l\leq x\leq u\right\}\ .$$

\vskip.5cm
\subsubsection*{\sc Convex integer maximization}
\label{cim}

We proceed to discuss the {\em maximization} of a convex function over
of the composite form $f(Wx)$, with $f:\Z^d\rightarrow\Z$ any
convex function and $W$ any integer $d\times n$ matrix.

\begin{figure}[hbt]
\hskip-1.1cm
\includegraphics[scale=0.58]{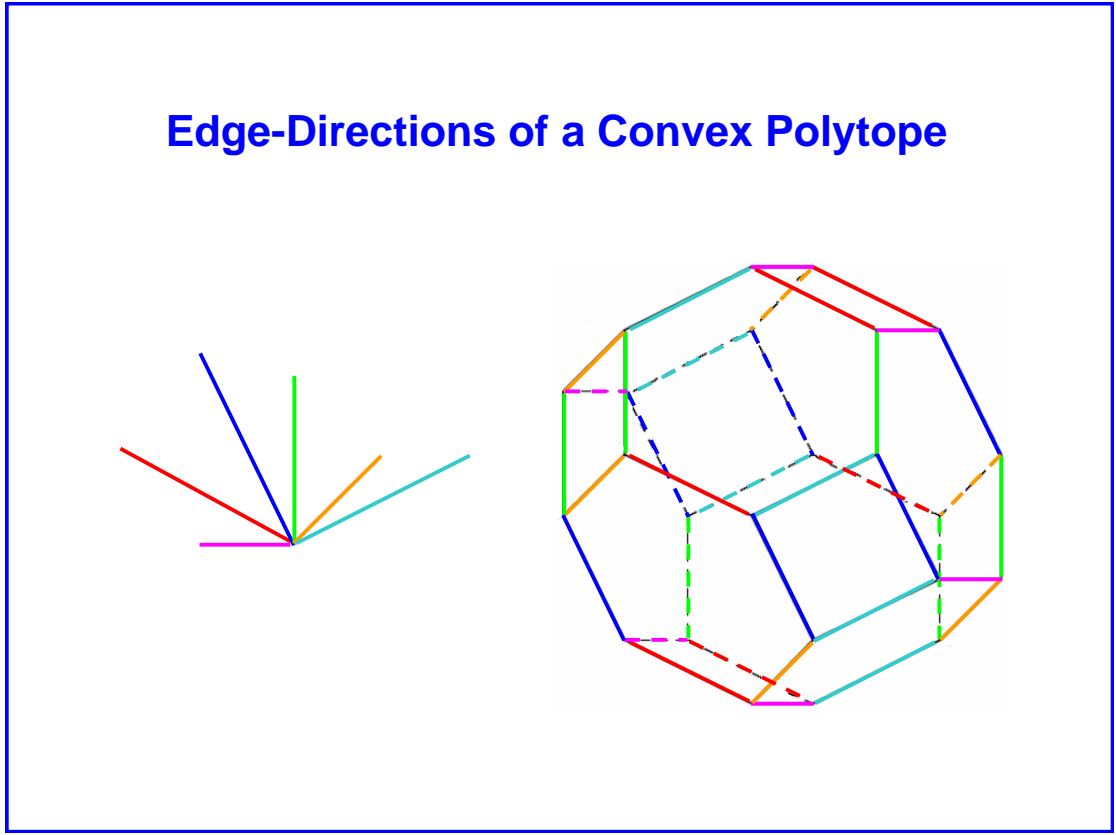}
\caption{Edge-Directions of a Convex Polytope}
\label{edge_directions_figure}
\end{figure}

We need a result of \cite{OR}. A {\em linear-optimization oracle}
for a set $S\subset\Z^n$ is one that,
given $w\in\Z^n$, solves the linear optimization problem $\max\{wx\,:\,x\in S\}$.
A {\em direction} of an edge ($1$-dimensional face) $e$ of a polyhedron $P$
is any nonzero scalar multiple of $u-v$ where $u,v$ are any two distinct
points in $e$. A set of {\em all edge-directions of $P$} is one
that contains some direction of each edge of $P$,
see Figure \ref{edge_directions_figure}.
\bt{binary_convex_maximization}{\bf \cite{OR}}
For all fixed $d$ there is an algorithm that, given a finite set $S\subset\Z^n$
presented by linear-optimization oracle, integer $d\times n$ matrix $W$, set
$E\subset\Z^n$ of all edge-directions of $\conv(S)$, and convex
$f:\Z^d\rightarrow\R$ presented by comparison oracle, solves in time
polynomial in $\l\max\{\|x\|_\infty:x\in S\},W,E\r$, the convex problem
$$\max\, \{f(Wx)\ :\ x\in S\}\ .$$
\et
We now show that, fortunately enough, the Graver basis of a
matrix $A$ is a set of all edge-directions of the integer hull
related to the integer program defined by $A$.
\bl{GraverEdges}
For every integer $m\times n$ matrix $A$, $l,u\in\Z_{\infty}^n$, and
$b\in\Z^m$, the Graver basis $\G(A)$ is a set of all
edge-directions of $P_I:=\conv\{x\in\Z^n:Ax=b,\ l\leq x\leq u\}$.
\el
\bpr
Consider any edge $e$ of $P_I$ and pick two distinct integer points $x,y\in e$.
Then $g:=y-x$ is in $\L^*(A)$ and hence Lemma \ref{ConformalGraverSumWeakForm}
implies that $g=\sum_i h_i$ is a conformal sum
for suitable $h_i\in\G(A)$. We claim that $x+h_i\in P_I$ for all $i$.
Indeed, $h_i\in\G(A)$ implies $A(x+h_i)=Ax=b$,
and $l\leq x,x+g \leq u$ and $h_i\sqsubseteq g$ imply $l\leq x+h_i\leq u$.

Now let $w\in\Z^n$ be uniquely maximized over $P_I$ at the edge $e$.
Then $wh_i=w(x+h_i)-wx\leq 0$ for all $i$.
But $\sum wh_i=wg=wy-wx=0$, implying that in fact $wh_i=0$
and hence $x+h_i\in e$ for all $i$. This implies that $h_i$
is a direction of $e$ (in fact, all $h_i$ are the same
and $g$ is a multiple of some Graver basis element).
\epr

Using Theorems \ref{Graver1} and \ref{binary_convex_maximization} and
Lemma \ref{GraverEdges} we obtain the following theorem.

\bt{Graver4}{\bf \cite{DHORW}}
For every fixed $d$ there is an algorithm that, given integer
$m\times n$ matrix $A$, its Graver basis $\G(A)$, $l,u\in\Z_{\infty}^n$,
$b\in\Z^m$, integer $d\times n$ matrix $W$, and convex
function $f:\Z^d\rightarrow\R$ presented by a comparison oracle, solves in time
which is polynomial in $\l A,W,\G(A),l,u,b\r$, the convex integer maximization problem
$$\max\, \{f(Wx)\ :\ x\in\Z^n,\ Ax=b,\ l\leq x\leq u\}\ .$$
\et
\bpr
Let $S:=\{x\in\Z^n\,:\,Ax=b\,,\, l\leq x\leq u\}$. The algorithm of
Theorem \ref{Graver1} allows to simulate in polynomial time a linear-optimization
oracle for $S$. In particular, it allows to either conclude that $S$ is infinite and
stop or conclude that it is finite, in which case $\l\max\{\|x\|_\infty:x\in S\}\r$
is polynomial in $\l A,l,u,b\r$, and continue. By Lemma \ref{GraverEdges},
the given Graver basis is a set of all edge-directions of $\conv(S)=P_I$.
Hence the algorithm of Theorem \ref{binary_convex_maximization} can be applied,
and provides the polynomial time solution of the convex integer maximization program.
\epr

\subsection{N-Fold Integer Programming}
\label{nfip}

In this subsection we focus our attention on (nonlinear) $n$-fold integer
programming. In \S \ref{gbonfp} we study Graver bases of $n$-fold products
of integer bimatrices and show that they can be computed in polynomial time.
In \S \ref{nfipipt} we combine the results of \S \ref{gbanip} and \S \ref{gbonfp},
and prove our Theorems \ref{NFold1}--\ref{NFold5}, which establish the
polynomial time solvability of linear and nonlinear $n$-fold integer programming.

\subsubsection{Graver bases of n-fold products}
\label{gbonfp}

Let $A$ be a fixed integer $(r,s)\times t$ bimatrix with blocks $A_1$, $A_2$.
For each positive integer $n$ we index vectors in $\Z^{nt}$ as
$x=(x^1,\ldots,x^n)$ with each {\em brick} $x^k$ lying in $\Z^t$.
The {\em type} of vector $x$ is the number
$\type(x):=|\{k\,:\,x^k\neq 0\}|$ of nonzero bricks of $x$.

The following definition plays an important role in the sequel.
\bd{GraverComplexityDefinition}{\bf \cite{SS}}
The {\em Graver complexity} of an integer bimatrix $A$ is defined as
$$g(A)\quad :=\quad \inf\left\{g\in\Z_+\ :\ \type(x)\leq g\ \
\mbox{for all $x\in\G(\A)$ and all $n$}\right\}\quad .$$
\ed
We proceed to establish a result of \cite{SS} and its
extension in \cite{HSu} which show that, in fact, the
Graver complexity of every integer bimatrix $A$ is finite.

Consider $n$-fold products $\A$ of $A$. By definition of the $n$-fold product,
$A^{(n)}x=0$ if and only if $A_1\sum_{k=1}^nx^k=0$ and $A_2x^k=0$ for all $k$.
In particular, a necessary condition for $x$ to lie in $\L(\A)$,
and in particular in $\G(\A)$, is that $x^k\in\L(A_2)$
for all $k$. Call a vector $x=(x^1,\dots,x^n)$ {\em full} if, in fact, $x^k\in\L^*(A_2)$
for all $k$, in which case $\type(x)=n$, and {\em pure} if, moreover,
$x^k\in\G(A_2)$ for all $k$. Full vectors, and in particular pure
vectors, are natural candidates for lying in the Graver basis $\G(\A)$
of $\A$, and will indeed play an important role in its construction.

Consider any full vector $y=(y^1,\dots,y^m)$. By definition, each brick
of $y$ satisfies $y^i\in\L^*(A_2)$ and is therefore a conformal sum
$y^i=\sum_{j=1}^{k_i} x^{i,j}$ of some elements $x^{i,j}\in\G(A_2)$
for all $i,j$. Let $n:=k_1+\cdots+k_m\geq m$ and let $x$ be the pure vector
$$x\ =\ (x^1,\dots,x^n)\ :=\
(x^{1,1},\dots,x^{1,k_1},\dots,x^{m,1},\dots,x^{m,k_m})\ .$$
We call the pure vector $x$ an {\em expansion} of the full vector $y$,
and we call the full vector $y$ a {\em compression} of the pure vector $x$.
Note that $A_1\sum y^i=A_1\sum x^{i,j}$ and therefore
$y\in\L(A^{(m)})$ if and only if $x\in\L(\A)$.
Note also that each full $y$ may have many different expansions and
each pure $x$ may have many different compressions.
\bl{Expansion}
Consider any full $y=(y^1,\dots,y^m)$ and any
expansion $x=(x^1,\dots,x^n)$ of $y$.
If $y$ is in the Graver basis $\G(A^{(m)})$
then $x$ is in the Graver basis $\G(\A)$.
\el
\bpr
Let $x=(x^{1,1},\dots,x^{m,k_m})=(x^1,\dots,x^n)$ be an expansion
of $y=(y^1,\dots,y^m)$ with $y^i=\sum_{j=1}^{k_i} x^{i,j}$ for each $i$.
Suppose indirectly $y\in\G(A^{(m)})$ but $x\notin\G(\A)$.
Since $y\in\L^*(A^{(m)})$ we have $x\in\L^*(\A)$. Since $x\notin\G(\A)$,
there exists an element $g=(g^{1,1},\dots,g^{m,k_m})$ in $\G(\A)$ satisfying
$g\sqsubset x$. Let $h=(h^1,\dots,h^m)$ be the compression of $g$ defined by
$h^i:=\sum_{j=1}^{k_i} g^{i,j}$. Since $g\in\L^*(\A)$ we have $h\in\L^*(A^{(m)})$.
But $h\sqsubset y$, contradicting $y\in\G(A^{(m)})$.
This completes the proof.
\epr
\bl{GraverComplexity}
The Graver complexity $g(A)$ of every integer bimatrix $A$ is finite.
\el
\bpr
We need to bound the type of any element in the Graver basis of the $l$-fold
product of $A$ for any $l$. Suppose there is an element $z$ of type $m$
in some $\G(A^{(l)})$. Then its restriction $y=(y^1,\ldots,y^m)$ to its $m$
nonzero bricks is a full vector and is in the Graver basis $\G(A^{(m)})$. Let
$x=(x^1,\dots,x^n)$ be any expansion of $y$. Then $\type(z)=m\leq n=\type(x)$,
and by Lemma \ref{Expansion}, the pure vector $x$ is in $\G(\A)$.

Therefore, it suffices to bound the type of any pure element
in the Graver basis of the $n$-fold product of $A$ for any $n$.
Suppose $x=(x^1,\ldots,x^n)$ is a pure element in $\G(A^{(n)})$ for some $n$.
Let $\G(A_2)=\{g^1,\dots,g^p\}$ be the Graver basis of $A_2$
and let $G_2$ be the $t\times p$ matrix whose columns are the $g^i$.
Let $v\in\Z_+^p$ be the vector with $v_i:=|\{k:x^k=g^i\}|$ counting
the number of bricks of $x$ which are equal to $g^i$ for each $i$.
Then $\sum_{i=1}^p v_i=\type(x)=n$. Now, note that
$A_1G_2v=A_1\sum_{k=1}^n x^k=0$ and hence $v\in\L^*(A_1G_2)$. We claim that,
moreover, $v$ is in $\G(A_1G_2)$. Suppose indirectly it is not. Then there
is a ${\hat v}\in\G(A_1G_2)$ with ${\hat v}\sqsubset v$, and it is easy to
obtain a nonzero ${\hat x}\sqsubset x$ from $x$ by zeroing out
some bricks so that ${\hat v}_i=|\{k:{\hat x}^k=g^i\}|$
for all $i$. Then $A_1\sum_{k=1}^n {\hat x}^k=A_1G_2{\hat v}=0$
and hence ${\hat x}\in\L^*  (A^{(n)})$, contradicting $x\in\G(A^{(n)})$.

So the type of any pure vector, and hence the Graver
complexity of $A$, is at most the largest value $\sum_{i=1}^p v_i$
of any nonnegative vector $v$ in the Graver basis $\G(A_1G_2)$.
\epr
We proceed to establish the following theorem from \cite{DHOW} which asserts
that Graver bases of $n$-fold products can be computed in polynomial time.
An {\em $n$-lifting} of a vector $y=(y^1,\dots,y^m)$ consisting of $m$ bricks
is any vector $z=(z^1,\dots,z^n)$ consisting of $n$ bricks such that for some
$1\leq k_1<\cdots<k_m\leq n$ we have $z^{k_i}=y^i$ for $i=1,\dots,m$, and all
other bricks of $z$ are zero; in particular, $n\geq m$ and $\type(z)=\type(y)$.
\bt{GraverComputation}{\bf \cite{DHOW}}
For every fixed integer bimatrix $A$ there is an algorithm that,
given positive integer $n$, computes the Graver basis $\G(\A)$
of the $n$-fold product of $A$, in time which is polynomial in $n$.
In particular, the cardinality $|\G(\A)|$ and the binary-encoding length
$\l \G(\A)\r$ of the Graver basis of $\A$ are polynomial in $n$.
\et
\bpr
Let $g:=g(A)$ be the Graver complexity of $A$. Since $A$ is fixed, so is $g$.
Therefore, for every $n\leq g$, the Graver basis $\G(\A)$, and in particular,
the Graver basis $\G(A^{(g)})$ of the $g$-fold product of $A$,
can be computed in constant time.

Now, consider any $n>g$. We claim that $\G(\A)$ satisfies
$$\G(\A)\ =\
\left\{z\ :\ z\ \mbox{is an $n$-lifting of some}\ y\in \G(A^{(g)})\right\}\ .$$
Consider any $n$-lifting $z$ of any $y\in\G(A^{(g)})$.
Suppose indirectly $z\notin\G(\A)$. Then there exists $z'\in\G(\A)$ with
$z'\sqsubset z$. But then $z'$ is the $n$-lifting of some
$y'\in\L^*(A^{(g)})$ with $y'\sqsubset y$,
contradicting $y\in\G(A^{(g)})$. So $z\in\G(\A)$.

Conversely, consider any $z\in\G(\A)$. Then $\type(z)\leq g$ and hence $z$
is the $n$-lifting of some $y\in\L^*(A^{(g)})$. Suppose indirectly
$y\notin\G(A^{(g)})$. Then there exists $y'\in\G(A^{(g)})$
with $y'\sqsubset y$. But then the $n$-lifting $z'$
of $y'$ satisfies $z'\in\L^*(\A)$ with $z'\sqsubset z$,
contradicting $z\in\G(\A)$. So $y\in\G(A^{(g)})$.

Now, the number of $n$-liftings of each $y\in\G(A^{(g)})$ is at
most ${n\choose g}$, and hence
$$|\G(\A)|\ \leq\ {n\choose g}|\G(A^{(g)})|\ =\ O(n^g)\ .$$
So the set of all $n$-liftings of vectors in $\G(A^{(g)})$ and hence
the Graver basis $\G(\A)$ of the $n$-fold product can be computed
in time polynomial in $n$ as claimed.
\epr

\subsubsection{N-fold integer programming in polynomial time}
\label{nfipipt}

Combining Theorem \ref{GraverComputation} and the results
of \S \ref{gbanip} we now obtain Theorems \ref{NFold1}--\ref{NFold4}.

\vskip.2cm\noindent{\bf Theorem \ref{NFold1} \cite{DHOW}}
{\em
For each fixed integer $(r,s)\times t$ bimatrix $A$, there is
an algorithm that, given positive integer $n$, $l,u\in\Z_{\infty}^{nt}$,
$b\in\Z^{r+ns}$, and $w\in\Z^{nt}$, solves in time which is polynomial in $n$ and
$\l l,u,b,w\r$, the following linear $n$-fold integer program,
$$\min\,\left\{wx\ :\ x\in\Z^{nt}\,,\ A^{(n)} x=b\,,\ l\leq x\leq u\right\}\ .$$
}
\vskip.2cm
\bpr
Compute the Graver basis $\G(\A)$ using the algorithm of
Theorem \ref{GraverComputation}. Now apply the algorithm of
Theorem \ref{Graver1} with this Graver basis and solve the problem.
\epr
\vskip.2cm\noindent{\bf Theorem \ref{NFold2} \cite{HOW1}}
{\em
For each fixed integer $(r,s)\times t$ bimatrix $A$, there is an algorithm
that, given $n$, $l,u\in\Z_{\infty}^{nt}$, $b\in\Z^{r+ns}$,
and separable convex $f:\Z^{nt}\rightarrow\Z$
presented by a comparison oracle, solves in time polynomial in $n$ and
$\l l,u,b,{\hat f}\r$, the program
\begin{equation*}
\min\left\{f(x)\ :\ x\in\Z^{nt}\,,\ A^{(n)} x=b\,,\ l\leq x\leq u\right\}\ .
\end{equation*}
}
\vskip.2cm
\bpr
Compute the Graver basis $\G(\A)$ using the algorithm of
Theorem \ref{GraverComputation}. Now apply the algorithm of
Theorem \ref{Graver2} with this Graver basis and solve the problem.
\epr
\vskip.2cm\noindent{\bf Theorem \ref{NFold3} \cite{HOW1}}
{\em
For each fixed integer $(r,s)\times t$ bimatrix $A$, there is an algorithm
that, given positive integers $n$ and $p$, $l,u\in\Z_{\infty}^{nt}$, $b\in\Z^{r+ns}$,
and ${\hat x}\in\Z^{nt}$, solves in time polynomial in $n$, $p$,
and $\l l,u,b,{\hat x}\r$, the following distance minimization program,
\begin{equation}\label{distance_equation_third}
\min\,\{\|x-{\hat x}\|_p\ :\ x\in\Z^{nt},\ \A x=b,\ l\leq x\leq u\}\ .
\end{equation}
For $p=\infty$ the problem (\ref{distance_equation_third})
can be solved in time polynomial in $n$ and $\l l,u,b,{\hat x}\r$.
}
\vskip.2cm
\bpr
Compute the Graver basis $\G(\A)$ using the algorithm of
Theorem \ref{GraverComputation}. Now apply the algorithm of
Theorem \ref{Graver3} with this Graver basis and solve the problem.
\epr
\vskip.2cm\noindent{\bf Theorem \ref{NFold4} \cite{DHORW}}
{\em
For each fixed $d$ and $(r,s)\times t$ integer bimatrix $A$,
there is an algorithm that, given $n$, bounds
$l,u\in\Z_{\infty}^{nt}$, integer $d\times {nt}$ matrix $W$,
$b\in\Z^{r+ns}$, and convex function $f:\Z^d\rightarrow\R$ presented
by a comparison oracle, solves in time polynomial in $n$ and $\l W,l,u,b \r$,
the convex $n$-fold integer maximization program
$$\max\{f(Wx)\ :\ x\in\Z^{nt}\,,\ A^{(n)}x=b\,,\ l\leq x\leq u\}\ .$$
}
\vskip.2cm
\bpr
Compute the Graver basis $\G(\A)$ using the algorithm of
Theorem \ref{GraverComputation}. Now apply the algorithm of
Theorem \ref{Graver4} with this Graver basis and solve the problem.
\epr

\subsubsection{Weighted separable convex integer minimization}
\label{wscim}

We proceed to establish Theorem \ref{NFold5} which is a broad
extension of Theorem \ref{NFold2} that allows the objective
function to include a composite term of the form $f(Wx)$, where
$f:\Z^d\rightarrow\Z$ is a separable convex function and $W$ is
an integer matrix with $d$ rows, and to incorporate inequalities on $Wx$.
We begin with two lemmas. As before, ${\hat f},{\hat g}$ denote
the maximum values of $|f(Wx)|,|g(x)|$ over the feasible set.
\bl{Graver5}
There is an algorithm that, given an integer $m\times n$
matrix $A$, an integer $d\times n$ matrix $W$,
$l,u\in\Z_{\infty}^n$, ${\hat l},{\hat u}\in\Z_{\infty}^d$, $b\in\Z^m$,
the Graver basis $\G(B)$ of
$$B\ :=\ \left(
\begin{array}{cc}
  A  &  0  \\
  W  & I  \\
\end{array}
\right)\ ,$$
and separable convex functions $f:\Z^d\rightarrow\Z$,
$g:\Z^n\rightarrow\Z$ presented by evaluation oracles, solves in
time polynomial in $\l A,W,\G(B),l,u,{\hat l},{\hat u},b,{\hat f},{\hat g}\r$,
the problem
\begin{equation}\label{StrongSeparableConvexEquation}
\min\{f(Wx)+g(x)\ :\ x\in\Z^n\,,\ Ax=b
\,,\ {\hat l}\leq Wx\leq{\hat u}\,,\ l\leq x\leq u\}\ .
\end{equation}
\el
\bpr
Define $h:\Z^{n+d}\rightarrow\Z$ by $h(x,y):=f(-y)+g(x)$ for all
$x\in\Z^n$ and $y\in\Z^d$. Clearly, $h$ is separable convex since $f,g$ are.
Now, problem (\ref{StrongSeparableConvexEquation}) can be rewritten as
$$\min\{h(x,y):(x,y)\in\Z^{n+d},\
\left(
\begin{array}{cc}
  A  &  0  \\
  W  &  I  \\
\end{array}
\right)
\left(
\begin{array}{c}
  x  \\
  y  \\
\end{array}
\right)=
\left(
\begin{array}{c}
  b  \\
  0  \\
\end{array}
\right),\
l\leq x\leq u,-{\hat u}\leq y\leq-{\hat l}\}\ ,$$
and the statement follows at once by applying Theorem
\ref{Graver2} to this problem.
\epr
\bl{ExtendedGraverComputation}
For every fixed integer $(r,s)\times t$ bimatrix $A$ and $(p,q)\times t$ bimatrix $W$,
there is an algorithm that, given any positive integer $n$, computes in time
polynomial in $n$, the Graver basis $\G(B)$ of the following
$(r+ns+p+nq)\times(nt+p+nq)$ matrix,
$$B\ :=\ \left(
\begin{array}{cc}
  \A  &  0  \\
  W^{(n)}  & I  \\
\end{array}
\right)\ .$$
\el
\bpr
Let $D$ be the $(r+p,s+q)\times (t+p+q)$ bimatrix whose blocks are defined by
$$
D_1\ :=\
\left(
\begin{array}{ccc}
A_1 & 0   & 0 \\
W_1 & I_p & 0 \\
\end{array}
\right)\ ,\quad
D_2\ :=\
\left(
\begin{array}{ccc}
A_2 & 0 & 0   \\
W_2 & 0 & I_q \\
\end{array}
\right)\quad.
$$
Apply the algorithm of Theorem \ref{GraverComputation} and
compute in polynomial time the Graver basis $\G(D^{(n)})$ of
the $n$-fold product of $D$, which is the following matrix:
$$
D^{(n)}\ =\
{\small\left(
\begin{array}{ccc|ccc|c|ccc}
A_1 & 0   & 0 & A_1 & 0   & 0 & \cdots & A_1 & 0   & 0 \\
W_1 & I_p & 0 & W_1 & I_p & 0 & \cdots & W_1 & I_p & 0 \\
\hline
A_2 & 0 & 0   & 0   & 0   & 0 & \cdots & 0   & 0   & 0 \\
W_2 & 0 & I_q & 0   & 0   & 0 & \cdots & 0   & 0   & 0 \\
\hline
0   & 0 & 0   & A_2 & 0 & 0   & \cdots & 0   & 0   & 0 \\
0   & 0 & 0   & W_2 & 0 & I_q & \cdots & 0   & 0   & 0 \\
\hline
\vdots & \vdots & \vdots & \vdots & \vdots & \vdots & \ddots & \vdots & \vdots & \vdots\\
\hline
0   & 0 & 0   & 0   & 0 & 0   & \cdots & A_2 & 0   & 0 \\
0   & 0 & 0   & 0   & 0 & 0   & \cdots & W_2 & 0   & I_q \\
\end{array}
\right)}\ .
$$
Suitable row and column permutations applied to $D^{(n)}$ give the following matrix:
$$
C\ :=\
{\small\left(
\begin{array}{cccc|cccc|cccc}
A_1 & A_1 & \cdots & A_1 & 0   & 0   & \cdots & 0    & 0   & 0   & \cdots & 0   \\
A_2 & 0   & \cdots & 0   & 0   & 0   & \cdots & 0    & 0   & 0   & \cdots & 0   \\
0   & A_2 & \cdots & 0   & 0   & 0   & \cdots & 0    & 0   & 0   & \cdots & 0   \\
\vdots & \vdots & \ddots & \vdots & \vdots & \vdots  &
\ddots & \vdots & \vdots & \vdots & \ddots & \vdots \\
0   & 0   & \cdots & A_2 & 0   & 0   & \cdots & 0    & 0   & 0   & \cdots & 0   \\
\hline
W_1 & W_1 & \cdots & W_1 & I_p & I_p & \cdots & I_p  & 0   & 0   & \cdots & 0   \\
W_2 & 0   & \cdots & 0   & 0   & 0   & \cdots & 0    & I_q & 0   & \cdots & 0   \\
0   & W_2 & \cdots & 0   & 0   & 0   & \cdots & 0    & 0   & I_q & \cdots & 0   \\
\vdots & \vdots & \ddots & \vdots & \vdots & \vdots &
\ddots & \vdots & \vdots & \vdots & \ddots & \vdots \\
0   & 0   & \cdots & W_2 & 0   & 0   & \cdots & 0    & 0   & 0   & \cdots & I_q \\
\end{array}
\right)}\ .
$$
Obtain the Graver basis $\G(C)$ in polynomial time from $\G(D^{(n)})$
by permuting the entries of each element of the latter by the
permutation of the columns of $\G(D^{(n)})$ that is used to get $C$
(the permutation of the rows does not affect the Graver basis).

Now, note that the matrix $B$ can be obtained from $C$ by dropping all
but the first $p$ columns in the second block. Consider any element
in $\G(C)$, indexed, according to the block structure, as
$(x^1,x^2,\dots,x^n,y^1,y^2,\dots,y^n,z^1,z^2,\dots,z^n)$.
Clearly, if $y^k=0$ for $k=2,\dots,n$ then the restriction
$(x^1,x^2,\dots,x^n,y^1,z^1,z^2,\dots,z^n)$ of this element is in the
Graver basis of $B$. On the other hand, if
$(x^1,x^2,\dots,x^n,y^1,z^1,z^2,\dots,z^n)$ is any element in $\G(B)$
then its extension $(x^1,x^2,\dots,x^n,y^1,0,\dots,0,z^1,z^2,\dots,z^n)$
is clearly in $\G(C)$. So the Graver basis of $B$ can be obtained
in polynomial time by
$$\G(B)\, :=\,\left\{(x^1,\dots,x^n,y^1,z^1,\dots,z^n)\, :\,
(x^1,\dots,x^n,y^1,0,\dots,0,z^1,\dots,z^n)\in\G(C)\right\}\, .
$$
This completes the proof.
\epr

\vskip.2cm\noindent{\bf Theorem \ref{NFold5} \cite{HOW2}}
{\em
For each fixed integer $(r,s)\times t$ bimatrix $A$ and integer $(p,q)\times t$
bimatrix $W$, there is an algorithm that, given $n$,
$l,u\in\Z_{\infty}^{nt}$, ${\hat l},{\hat u}\in\Z_{\infty}^{p+nq}$, $b\in\Z^{r+ns}$,
and separable convex
functions $f:\Z^{p+nq}\rightarrow\Z$, $g:\Z^{nt}\rightarrow\Z$ presented
by evaluation oracles, solves in time
polynomial in $n$ and $\l l,u,{\hat l},{\hat u},b,{\hat f},{\hat g}\r$,
the generalized program
\begin{equation*}\label{NFoldStrongSeparableConvexEquation}
\min\left\{f(W^{(n)}x)+g(x)\ :\ x\in\Z^{nt}\,,\ A^{(n)}x=b\,,\
{\hat l}\leq W^{(n)}x\leq{\hat u}\,,\ l\leq x\leq u\right\}\ .
\end{equation*}
}

\vskip.2cm
\bpr
Use the algorithm of Lemma \ref{ExtendedGraverComputation}
to compute the Graver basis $\G(B)$ of
$$B\ :=\ \left(
\begin{array}{cc}
  \A  &  0  \\
  W^{(n)}  & I  \\
\end{array}
\right)\ .$$
Now apply the algorithm of Lemma \ref{Graver5}
and solve the nonlinear integer program.
\epr

\section{Discussion}
\label{d}

We conclude with a short discussion of the universality of
$n$-fold integer programming and the Graver complexity of (directed) graphs,
a new important invariant which controls the complexity of our
multiway table and multicommodity flow applications.

\subsection{Universality of N-Fold Integer Programming}
\label{uonfip}

Let us introduce the following notation. For an integer $s\times t$
matrix $D$, let $\boxminus D$ denote the $(t,s)\times t$ bimatrix
whose first block is the $t\times t$ identity matrix and whose second block is $D$.
Consider the following special form of the $n$-fold product, defined
for a matrix $D$, by $D^{[n]}:=\left(\boxminus D\right)^{(n)}$.
We consider such $m$-fold products of the $1\times 3$ matrix ${\bf 1}_3:=[1,1,1]$.
Note that ${\bf 1}_3^{[m]}$ is precisely the $(3+m)\times 3m$ incidence
matrix of the complete bipartite graph $K_{3,m}$.
For instance, for $m=3$, it is the matrix
\begin{equation*}\label{matrix}
{\bf 1}_3^{[3]}\ =\
\left(
\begin{array}{ccccccccc}
  1 & 0 & 0 & 1 & 0 & 0 & 1 & 0 & 0 \\
  0 & 1 & 0 & 0 & 1 & 0 & 0 & 1 & 0 \\
  0 & 0 & 1 & 0 & 0 & 1 & 0 & 0 & 1 \\
  1 & 1 & 1 & 0 & 0 & 0 & 0 & 0 & 0 \\
  0 & 0 & 0 & 1 & 1 & 1 & 0 & 0 & 0 \\
  0 & 0 & 0 & 0 & 0 & 0 & 1 & 1 & 1 \\
\end{array}
\right)\ .
\end{equation*}
We can now rewrite Theorem \ref{TableUniversality}
in the following compact and elegant form.

\vskip.2cm\noindent
{\bf The Universality Theorem \cite{DO2}}
{\em Every rational polytope $\{y\in\R_+^d\,:\,Ay=b\}$ stands in polynomial
time computable integer preserving bijection with some polytope
\begin{equation}\label{universality_equation}
\left\{x\in\R_+^{3mn}\ :\ {\bf 1}_3^{[m][n]}x=a\right\}\ .
\end{equation}}

\vskip.2cm\noindent
The bijection constructed by the algorithm of this theorem is,
moreover, a simple projection from $\R^{3mn}$ to $\R^d$
that erases all but some $d$ coordinates (see \cite{DO2}).
For $i=1,\dots,d$ let $x_{\sigma(i)}$ be the coordinate of $x$
that is mapped to $y_i$ under this projection. Then any linear or
nonlinear integer program $\min\{f(y)\,:\,y\in\Z_+^d,\,Ay=b\}$
can be lifted in polynomial time to the following integer program
over a simple $\{0,1\}$-valued matrix ${\bf 1}_3^{[m][n]}$
which is completely determined by two parameters $m$ and $n$ only,
\begin{equation}\label{universality_optimization_equation}
\min\left\{f\left(x_{\sigma(1)},\dots,x_{\sigma(d)}\right)\ :\
x\in\Z_+^{3mn}\,,\ {\bf 1}_3^{[m][n]}x=a\right\}\ .
\end{equation}
This also shows the universality of $n$-fold integer programming:
every linear or nonlinear integer program is equivalent to an
$n$-fold integer program over some bimatrix $\boxminus{\bf 1}_3^{[m]}$
which is completely determined by a single parameter $m$.

Moreover, for every fixed $m$, program (\ref{universality_optimization_equation})
can be solved in polynomial time for linear forms and broad classes of
convex and concave functions by Theorems \ref{NFold1}--\ref{NFold5}.

\subsection{Graver Complexity of Graphs and Digraphs}
\label{gcogad}

The significance of the following new (di)-graph invariant will be explained below.
\bd{GraphGraverComplexityDefinition}{\bf \cite{BO}}
The {\em Graver complexity} of a graph or a digraph $G$
is the Graver complexity $g(G):=g(\boxminus D)$
of the bimatrix $\boxminus D$ with $D$ the incidence matrix of $G$.
\ed

One major task done by our algorithms for linear and nonlinear $n$-fold
integer programming over a bimatrix $A$ is the construction of the
Graver basis $\G(\A)$ in time $O\left(n^{g(A)}\right)$ with $g(A)$
the Graver complexity of $A$ (see proof of Theorem \ref{GraverComputation}).

Since the bimatrix underlying the universal $n$-fold integer
program (\ref{universality_optimization_equation}) is precisely
$\boxminus D$ with $D={\bf 1}_3^{[m]}$ the incidence matrix of $K_{3,m}$,
it follows that the complexity of computing the relevant Graver bases
for this program for fixed $m$ and variable $n$ is
$O\left(n^{g(K_{3,m})}\right)$ where $g(K_{3,m})$ is the
Graver complexity of $K_{3,m}$ as just defined.

Turning to the many-commodity transshipment problem over a digraph $G$
discussed in \S \ref{tmctp}, the bimatrix underlying the $n$-fold
integer program (\ref{many-commodity_equation})
in the proof of Theorem \ref{Transshipment} is precisely
$\boxminus D$ with $D$ the incidence matrix of $G$, and so
it follows that the complexity of computing the relevant Graver bases
for this program is $O\left(n^{g(G)}\right)$ where $g(G)$ is the
Graver complexity of the digraph $G$ as just defined.

So the Graver complexity of a (di)-graph controls the
complexity of computing the Graver bases of the relevant
$n$-fold integer programs, and hence its significance.

Unfortunately, our present understanding of the Graver complexity of (di)-graphs
is very limited and much more study is required. Very little is known
even for the complete bipartite graphs $K_{3,m}$: while $g(K_{3,3})=9$,
already $g(K_{3,4})$ is unknown. See \cite{BO} for more details and a
lower bound on $g(K_{3,m})$ which is exponential in $m$.

\section*{Acknowledgements}

I thank Jon Lee and Sven Leyffer for inviting me to write this
article. I am indebted to Jesus De Loera, Raymond Hemmecke,
Uriel Rothblum and Robert Weismantel for their collaboration in
developing the theory of $n$-fold integer programming, and to
Raymond Hemmecke for his invaluable suggestions.
The article was written mostly while I was visiting and
delivering the Nachdiplom Lectures at ETH Z\"urich during Spring 2009.
I thank the following colleagues at ETH for useful feedback:
David Adjiashvili, Jan Foniok, Martin Fuchsberger, Komei Fukuda,
Dan Hefetz, Hans-Rudolf K\"unsch, Hans-Jakob L\"uthi, and
Philipp Zumstein. I also thank Dorit Hochbaum, Peter Malkin,
and a referee for useful remarks.

\vskip0.5cm
\noindent {\small Shmuel Onn}\newline
\emph{Technion - Israel Institute of Technology, 32000 Haifa, Israel}\newline
and\newline
\emph{ETH Z\"urich, 8092 Z\"urich, Switzerland}\newline
\emph{onn{\small @}ie.technion.ac.il},\
\emph{http://ie.technion.ac.il/{\small $\sim$onn}}

\end{document}